\documentstyle[12pt]{article}
\newtheorem{ther}{Theorem}[section]
\newtheorem{teo}{Theorem}
\newtheorem{prob}{Problem}

\newtheorem{pr}{Proposition}[section]
\newtheorem{cor}{Corollary}[section]
\newtheorem{lm}{Lemma}[section]
\newtheorem{rem}{Remark}[section]
\newtheorem{ex}{Example}[section]
\newtheorem{exa}{Example}
\newcommand{\End}{{\mbox{End}\,}}
\newcommand{\Hom}{{\mbox{Hom}\,}}
\newcommand{\Lie}{{\mbox{Lie}\,}}
\newcommand{\ind}{{\mbox{ind}\,}}
\newcommand{\diag}{{\mbox{diag}\,}}

\newcommand{\Spec}{{\mbox{Spec}\,}}
\newcommand{\id}{{\mbox{id}\,}}
\title{Invariants of mixed representations of quivers I}
\author{A.N.Zubkov \\
644099, Omsk-99, Tuhachevskogo embankment 14,\\
Omsk State Pedagogical University ,\\
chair of algebra e-mail: zubkov@iitam.omsk.net.ru}
\date{}
\begin{document}
\maketitle

\section*{Introduction}

The concept of a representation of a quiver was introduced in \cite{gab}.
If we consider all representations of a quiver of given dimension as an affine
variety provided with the action of its automorphism group then the points of
the corresponding categorical quotient can be parametrized by semisimple
representations. Moreover, this quotient is also an affine variety and its
coordinate algebra is generated by all polynomial invariants.
In the characteristic zero case invariants of representations of
quivers were first described
in \cite{prb1, prb2}. This result was applied to
investigate an etale local structure of categorical quotients of
quiver representation spaces \cite{prb1, prb2}.

The modular case was explored in
\cite{don1, zub4}. In \cite{don1} invariants of arbitrary
quiver were described over any infinite field. In \cite{zub4}
all defining
relations between them are described too. We note that the last result 
was proved
independently
in \cite{dom} for the characteristic zero case.
Finally, in \cite{dz2} the main results from \cite{prb1, prb2} concerning an 
etale
local structure of invariants of a quiver were extended to the case of
any algebraically closed field.

No doubt the next step should be
to generalize these statements for other classical
groups, specifically to the orthogonal and symplectic groups.
It is clear that one has to start with
the action of $O(n)$ or $Sp(n)$ on $m$-tuples of $n\times n$
matrices by simultaneous conjugation. Using the so-called transfer principle
\cite{gr} one can reduce this problem to
a 
representation of some quiver.
This representation is a new type of
representations of quivers
called {\it mixed} representations.

Recall some necessary definitions and notations
(see \cite{gab, don1, prb1,
prb2}).
A {\it quiver} is a
quadruple $Q=(V,A,i,t)$, where $V$ is a vertex set and $A$ is an arrow set of
$Q$. Let the maps $i,t:A\rightarrow V$ associate to each arrow $a\in A$ its
origin $i(a)\in V$ and its end $t(a)\in V$. We enumerate elements of the
vertex set as $V=\{1,\ldots ,n\}$.

We consider a collection of vector
spaces $E_1,\ldots ,E_n$ over an algebraically closed field $K$. Set $\dim
E_1=d_1,\ldots ,\dim E_n=d_n$. Denote by ${\bf d}$ the vector $(d_1,\ldots
,d_n)$. This vector is called a {\it dimension vector}.
For two dimension vectors ${\bf d}(1),
{\bf d}(2)$ we write
${\bf d}(1)\geq{\bf d}(2)$ iff $\forall i\in V, d(1)_i\geq\d(2)_i$.

Denote by $GL({\bf d})$
the group $GL(E_1)\times \ldots \times GL(E_n)=GL(d_1)\times \ldots \times
GL(d_n)$. The representation space of a quiver $Q$ of dimension ${\bf d}$ is 
$R(Q, {\bf d})=\prod_{a\in A}\Hom_K(E_{i(a)},E_{t(a)})$. The group 
$GL({\bf d})$
acts on $R(Q,{\bf d})$ by the rule:

$$
(y_a)_{a\in A}^g=(g_{t(a)}y_ag_{h(a)}^{-1})_{a\in A},
g=(g_1,\ldots , g_n)\in GL({\bf d}),
$$

$$
(y_a)_{a\in A}\in R(Q, {\bf d}).
$$
 
\noindent
For example, if our quiver $Q$ has one vertex and $m$ loops which are
incident to
this vertex then the ${\bf d}=(d)$-representation space of this quiver is
isomorphic
to the space of $m$ $d\times d$-matrices with respect to the diagonal action
of the group $GL(d)$ by conjugation.

The coordinate ring of the affine variety $R(Q, {\bf d})$
is isomorphic to $K[y_{ij}(a)\mid 1\leq j\leq d_{i(a)}, 1\leq i\leq
d_{t(a)}, a\in A]$. For any $a\in A$ denote by $Y_{{\bf d}}(a)$ the general
matrix $(y_{ij}(a))_{1\leq j\leq d_{i(a)}, 1\leq i\leq d_{t(a)}}$. The
action of $GL({\bf d})$ on $R(Q, {\bf d})$ induces the action on the
coordinate ring by the rule $Y_{{\bf d}}(a)\mapsto g_{t(a)}^{-1}
Y_{{\bf d}}(a)g_{i(a)}, a\in A$. We omit the lower index ${\bf d}$ if it does
not
lead to confusion. For example, we write $Y(a)$ instead of
$Y_{{\bf d}}(a)$.

Let us partition the vertex set of the quiver $Q$ into
several disjoint subsets. To be precise, let  $V=
V_{ord}\bigsqcup(\bigsqcup_{q\in\Omega}V_q)$. The vertices from $V_{ord}$
are said to be {\it ordinary}. We require that 
all subsets $V_q$ have cardinality two,
that is for any $q\in\Omega$ $V_q=\{i_q, j_q\}$.

A dimension vector ${\bf d}$ is said to be {\it compatible} with this
partition of $V$ if for any
$q\in\Omega, d_{i_q}=d_{j_q}=d_q$.
From now on all dimension vectors
are compatible with some fixed partition $V=V_{ord}\bigsqcup(\bigsqcup_{q\in
\Omega}V_q)$  unless otherwise stated.

The next step is to replace all $E_{j_q}, q\in\Omega$ by their duals. To
indicate that some vertices correspond to the duals of vector spaces we
introduce a new dimension vector
${\bf t} =(t_1,\ldots , t_l)$, where $t_i=d_i$
iff we assign to $i$ the space $E_i$, otherwise $t_i=d_i^*$.
We call ${\bf d}$ the vector underlying ${\bf t}$. This notation
will be used throughout.

By definition, the ${\bf t}$-dimensional
representation space of the
quiver $Q$ is equal to the space $R(Q, {\bf t})=\prod_{a\in A}\Hom_K (W_{i(a)}
, W_{t(a)})$, where $W_i=E_i$ iff $t_i=d_i$, otherwise $W_i=E_i^*$.

The space $R(Q, {\bf t})$ is a $G=GL({\bf d})$-module
under the same action

$$(y_a)_{a\in A}^g=(g_{t(a)}y_ag_{i(a)}^{-1})_{a\in A},
g=(g_1,\ldots , g_l)\in G,$$

$$(y_a)_{a\in A}\in R(Q, {\bf t}).$$

\noindent If $\Omega=\emptyset$ then ${\bf t}={\bf d}$ and $R(Q, {\bf t})=
R(Q, {\bf d})$. Without loss of generality one can identify the coordinate 
algebras $K[R(Q, {\bf t})]$ and  $K[R(Q, {\bf d})]$.

Finally, replacing all subfactors $GL(E_{i_q})\times GL(E_{j_q})=GL(d_q)\times
GL(d_q)$ of the group $G=GL({\bf d})$  by
their diagonal subgroups we get a new group $H({\bf t})$.
The space $R(Q, {\bf t})$ with respect to the action of  the group
$H({\bf t})$ is called the {\it mixed} representation space of the quiver $Q$
of dimension ${\bf t}$ relative to the partition
$V=V_{ord}\bigsqcup(\bigsqcup_{q\in\Omega}V_q)$.

\begin{exa}
Let $V_{ord}=\emptyset, \Omega=\{q\}, i_q=1, j_q=2, A=\{a_1, \ldots , a_m\, 
b, c\}, i(a_k)=t(a_k)=i(b)=t(c), t(b)=i(c)=2$. The mixed representation
space of this quiver of dimension ${\bf t}=(d, d^*)$ can be identified with
$M(d)^m\times M(d)^2$, where the first $m$ $d\times d$ matrix coordinates
correspond to the loops $a_1, \ldots, a_m$ but the two last ones correspond
to the arrows $b, c$ respectively.
The group $GL({\bf t})=GL(d)$
acts on this mixed representation space by the rule

$$(A_1,\ldots , A_m, B, C)^g=(gA_1 g^{-1},\ldots , gA_m g^{-1}, (g^t)^{-1}
Bg^{-1}, gCg^t),$$
$$A_i, B, C\in M(d),
g\in GL(d), 1\leq i\leq m.
$$

\noindent This special case of mixed representations
of quivers first appeared in \cite{zub5} to compute the invariants of
orthogonal or symplectic groups acting diagonally by conjugations on 
several matrices.  
\end{exa}

We formulate the following.

\begin{prob}
What are the generators and the defining relations between them for the ring
$J(Q, {\bf t})=K[R(Q, {\bf t})]^{H({\bf t})}$?
\end{prob}

The principal aim of this article is to answer the first part of this
question as well as to prepare some necessary facts to answer  
the second part in the next article.

To formulate the main result of this article we need some additional 
definition.
Let us define a {\it doubled} quiver $Q^{(d)}$. The vertex set
$V^{(d)}$ of this
quiver is equal to $V\bigsqcup V^*_{ord}$, where $V^*_{ord}=
\{i^*\mid i\in V_{ord}\}$.
Respectively,
the arrow set $A^{(d)}$ of $Q^{(d)}$
is equal to $A\bigsqcup\overline{A}$, where $\overline{A}=\{\bar{a}\mid
a\in A\}$.
Further, if $i(a), t(a)\in V_{ord}$ then $i(\bar{a})=t(a)^*,
t(\bar{a})=i(a)^*$ but if $i(a)$ or $t(a)$ lies in some $V_q, q\in\Omega,$ then

$$i(\bar{a})=\left\{\begin{array}{c}
j_q, t(a)=i_q,\\
i_q, t(a)=j_q
\end{array} \right.
$$

and symmetrically

$$t(\bar{a})=\left\{\begin{array}{c}
j_q, i(a)=i_q,\\
i_q, i(a)=j_q
\end{array} \right.
.$$

\noindent Finally, for any $a\in A^{(d)}$ we suppose
$Z(a)=Y(a)$ if $a\in
A$ otherwise $a=\bar{b}, b\in A$ and $Z(a)=Y(b)^t$, where $Y(b)^t$ is
transpose of $Y(b)$.

A product $Z(a_m)\ldots Z(a_1)$ is said to be {\it admissible} if $a_m,
\ldots , a_1$ is a closed path in $Q^{(d)}$, that is if $t(a_i)=i(a_{i+1}),
i=1, \ldots, m-1$ and $i(a_1)=t(a_m)$.
A pair $Z(a)Z(b)$ is said to be {\it linked} if $t(b)=i(a)$. It is clear that
$Z(a_m)\ldots Z(a_1)$ is admissible iff all pairs $Z(a_{i+1})Z(a_i), i=1,
\ldots , m-1,$ and $Z(a_1)Z(a_m)$ are linked.

Using the theory of modules with good filtration as well as some reductions 
developed in \cite{don2, zub5} we prove

\begin{teo}
The algebra $J(Q, {\bf t})$ is generated by the elements $\sigma_j(Z(a_r)\ldots
Z(a_1))$, where $1\leq j\leq\max\limits_{1\leq i\leq n}\{ d_i\}$, $a_r,
\ldots, a_1$ is a closed path in the double quiver $Q^{(d)}$ and $\sigma_j$ is
$j$-th coefficient of characteristic polynomial.
\end{teo}

One can define more general {\it supermixed} representations of quivers 
involving as
special cases mixed representations of quivers and orthogonal (symplectic)
representations of symmetric quivers introduced in \cite{dw3}. To be precise,
let $R(Q, {\bf t})$ be the mixed representation space of a quiver $Q$ of
dimension ${\bf t}=(t_1,\ldots , t_l)$ with respect to some partition of
$V$, say $V=V_{ord}\bigsqcup (\bigsqcup_{q\in\Omega}V_q)$ as above.
By definition ${\bf t}$ is compatible with this
partition.

Replace some factors of the group $H=H({\bf t})=(\prod_{i\in V_{ord}}
GL(d_i))\times (\prod_{q\in\Omega} GL(d_q))$
by orthogonal or symplectic subgroups requiring additionally that the
characteristic of the ground field is odd if at least one factor
is replaced by an orthogonal group. Denote the subgroup of $H$ obtained in this
way as $G=(\prod_{i\in V_{ord}}G_i)\times (\prod_{q\in\Omega} G_q)$, where
each factor $G_i (G_q)$ is either general linear group, orthogonal group or
symplectic group of given dimension.

Next, let us extract among all components $\Hom_K(W_{h(a)}, W_{t(a)}),
a\in A$,  those having property $i(a), t(a)\in V_q, q\in\Omega$.
Let $i(a)=i, t(a)=j$.
We have three cases: $G_q=GL(d_q)$, $G_q=O(d_q)$ or $G_q=Sp(d_q)$.

Let us consider the first case $G_q=GL(d_q)$. Let $i=j_q, j=i_q$
or $i=i_q, j=j_q$, that is $t_i=d_q^*, t_j=d_q$ or $t_i=d_q, t_j=d_q^*$.
Identifying $\Hom_K(W_i, W_j)$ with $M(d_q)$
one can replace this space by its subspaces of symmetric or skew-symmetric
matrices. In notations of \cite{dw3} these subspaces can be identified
with $S^2(V) (S^2(V^*))$ or $\Lambda^2(V) (\Lambda^2(V^*))$ respectively in
obvious way as a $GL(d_q)$-modules, where $V=E_{i_q}=E_{j_q}$.

In two remaining cases it does not matter if $(t_i, t_j)$ coincide with
$(d_q^*, d_q)$ or with $(d_q, d_q^*)$. Indeed, $V\cong V^*$ as a $O(V)$ or
$Sp(V)$-module.
If $G_q=O(d_q)$ then one can replace the space $\Hom_K(W_i, W_j)=M(d_q)$
by its subspaces of symmetric or skew-symmetric matrices again.

In the case $G_q=Sp(d_q)$ one can replace the space
$\Hom_K(W_i, W_j)=M(d_q)$ by its subspaces $\Lie(Sp(d_q))=\{A\in M(d_q)\mid
AJ \ \mbox{is a symmetric matrix}\}$ or $\{A\in M(d_q)\mid
AJ \ \mbox{is a skew-symmetric matrix}\}$, where $J=J_{d_q}$ is a $d_q\times
d_q$ skew-symmetric matrice of the bilinear form defining the group $Sp(d_q)$.

Denote a subspace of $R(Q, {\bf t})$ obtained with the help of some
replacements described above by $S$.
A pair $(S, G)$ is said to be a {\it supermixed} representation
space of the quiver $Q$ with respect to the induced action of the group $G$.

\begin{exa}
The space of $m$ $d\times d$ matrices with respect to the diagonal action of
$O(d)$ or $Sp(d)$ by conjugatyions is a supermixed representation space of 
the quiver $Q$ with one vertex and $m$ loops incident to this vertex. 
\end{exa}

The invariants of the supermixed representation space from Example 2 can be
obtained by specialization of invariants of mixed representation space from
Example 1 \cite{zub5}. This case is typical. In fact, we prove

\begin{teo}
Let $(S, G)$ be a supermixed representation space of a quiver $Q$.
There exists a quiver $Q'$ such that the algebra 
$K[S]^G$ is an epimorphic image of
$K[R(Q', {\bf t})]^{H({\bf t})}$ for some dimension vector ${\bf t}$. 
\end{teo}

The mixed or supermixed representations of quivers naturally arose from
the actions of $O(n)$ or $Sp(n)$ on several $n\times n$
matrices by simultaneous conjugation. If we replace $O(n)$ by its subgroup
$SO(n)$ then we will have to investigate semi-invariants
of mixed representations of quivers. In other words,
one can set the problems to find the generators and defining
relations between them for semi-invariants of mixed representations of
quivers. 

The problem to describe semi-invariants of ordinary representations of quivers
was very popular
during the last 20 years starting with the remarkable Kac's article \cite{ka}.
Important results were obtained in \cite{sco, sco2}.
There is also an extensive literature on semi-invariants of
Dynkin and Euclidean (or extended Dynkin) quivers, see
\cite{abe}, \cite{rin}, \cite{koi1}, \cite{koi2},
\cite{hh}, \cite{scwe}, \cite{skw}.
The complete descriptions of semi-invariants for
an arbitrary quiver were obtained in \cite{dw1, dw2} and \cite{dz}. In the
characteristic zero case the similar result was proved in \cite{sv}.
I believe that the method of this article will also make possible to describe
semi-invariants of mixed or supermixed representations of quivers.   
 
\section{Preliminaries}

\subsection{Induced modules and good filtrations}

Let $G$ be an algebraic group, $H$ a closed subgroup of $G$, and $A$ a rational
$H$-module. Then $K[G]\otimes A$ is naturally a rational $G\times H$-module
with respect to the action $(g, h)\cdot f\otimes a=f^{(g, h)}\otimes ha$, where
$g\in G, h\in H, f\in K[G], a\in A$ and $f^{(g, h)}(x)=f(g^{-1}xh)$. The set of
$H$-fixed points is a rational $G$-submodule, called the induced module
$\ind^G_H A=(K[G]\otimes A)^H$ \cite{gr}.

\begin{pr}(\cite{gr}, Theorem 9.1) If $X$ is an affine $G$-variety, that is $G$
acts rationally on $X$, then the invariant algebra $K[X]^H$ is isomorphic to
$(K[X]\otimes k[G/H])^G$, where $G$ acts on $K[G/H]$ by left translation.
The isomorphism is given by $a\otimes f\mapsto af(eH)$.
\end{pr}

Let $G$ be a reductive group. Fix some maximal torus of the group $G$, say $T$,
and a Borel subgroup $B$ containing $T$. The group $B$ has a semi-direct
decomposition $B=T\bowtie U$, where $U$ is a maximal unipotent subgroup of the
group $B$.
Denote by $X(T)$ the character group of
the torus $T$ and by $X(T)^+$ the dominant weight subset of $X(T)$ corresponding
to $B$.
If $\mu\in X(T)^+$ then denote by $\bigtriangledown(\mu)$ the induced module
$\ind_{B^-}^{G} K_{\mu}$, where $B^-$ is the opposite Borel subgroup and
$K_{\mu}$ is the one-dimensional $B^-$-module with respect to the action
$(tu)\circ x=\mu(t)x, t\in T, u\in U^-, x\in K_{\mu}$.

We say that a $G$-module $V$ has a {\it good filtration} (briefly GF) 
if there is some filtration with at most countable number of members

$$0\subseteq V_1\subseteq V_2\subseteq\ldots , \bigcup_{i=1}^{\infty}V_i=V$$

\noindent such that $\forall i\geq 1, V_i/V_{i-1}\cong\bigtriangledown(\mu_i)$.
Respectively, we say that a $G$-module $W$ has a
{\it Weyl filtration} (briefly WF)  if
there is some filtration with at most countable number of members

$$0\subseteq W_1\subseteq W_2\subseteq\ldots , \bigcup_{i=1}^{\infty}W_i=W$$

\noindent such that $\forall i\geq 1, W_i/W_{i-1}\cong\bigtriangleup(\mu_i)$, 
where
$\bigtriangleup(\mu)\cong\bigtriangledown(\mu^*)^*, \mu^*=-w_0(\mu)$ and
$w_0$ is the longest element of the Weyl group $W(G, T)=N_G(T)/T$.

It is clear that a finite-dimensional $G$-module $V$ has WF iff the dual
module $V^*$ has GF. A finite-dimensional module $V$
is called a {\it tilting} one if both $V$ and $V^*$ are with GF. In other 
words, $V$
has good and Weyl filtrations simultaneously.

We list some standard properties of modules having GF \cite{jan, don3,
don5, mat1}.

\begin{ther}
\begin{enumerate}
\item If

$$0\rightarrow V\rightarrow W\rightarrow S\rightarrow0$$

is a short exact sequence of $G$-modules and $V$ has GF, then the diagram

$$0\mapsto V^G\rightarrow W^G\rightarrow S^G\rightarrow 0$$

is exact.

\item If $W$ is a $G$-module with GF and $V$ is a submodule of $W$ with GF,
then the quotient $W/V$ is also a $G$-module with GF.

\item For given $G$-modules with GF their tensor product with respect to the
diagonal action of the group $G$ is also a module with GF.

\item If $V$ is a $G$-module with GF and $H$ is a Levi subgroup or the
commutator subgroup of $G$, then $V$ has GF as a $H$-module.
\end{enumerate}
\end{ther}

\subsection{Necessary facts of representation theory of products of
general linear groups}

Let $G=GL(k)$ and $T(k)=\{\diag(t_1,\ldots , t_k)\mid t_1,\ldots ,
t_k\in K^*\}$ is the standard torus of $G$. We fix the Borel subgroup $B(k)$
consisting of all upper triangular matrices.
It is clear that $B^{-}(k)$ consists of all lower triangular matrices.

Any character $\lambda\in X(T(k))$
can be regarded as a vector $(\lambda_1,\ldots ,\lambda_k)$ with
integer coordinates. By definition $\lambda(t)=
t_1^{\lambda_1}\ldots t_k^{\lambda_k}, t\in T(k)$.
It is known that
$\lambda\in X(T(k))^+$ iff $\lambda_1\geq\ldots\geq\lambda_k$ \cite{don2}.
If additionally $\lambda_k\geq 0$ then $(\lambda_1,\ldots ,\lambda_k)$ is 
called
an {\it ordered partition} and $\bigtriangledown(\lambda)$ is isomorphic to
so-called {\it Schur} module $L_{\tilde{\lambda}}(K^k)$ (see the next
subsection), where
$\tilde{\lambda}$ is the partition conjugated to $\lambda$. 
To be precise, if $\lambda_1=
\ldots =\lambda_{s_1} >\lambda_{s_1+1}=\ldots =\lambda_{s_1+s_2} >\ldots
>\lambda_{s_1+\ldots s_f +1}=\ldots=\lambda_k$ then $\tilde{\lambda}=(
k^{\lambda_k}, (s_1+\ldots +s_f)^{\lambda_{s_f}},\ldots , s_1^{\lambda_{s_1}})$,
where $l^k$ means $\underbrace{l,\ldots , l}_{k}$.

\begin{ex} If $\lambda=
(1^t, 0^{k-t})$ then $\bigtriangledown(\lambda)=
L_{\tilde{\lambda}}(K^k)=\Lambda^t(K^k)$ the $t$-th exterior power of the space
$K^k$. Moreover, $\Lambda^t(K^k)^*\cong\Lambda^{k-t}(K^k)\otimes\det^{-1}=
\bigtriangledown(0^{k-t}, -1^t)$ has GF. In particular, 
$\Lambda^t(K^k)$ is
a tilting $GL(k)$-module. Using Theorem 1.1(3) we obtain that all tensor
products of such modules are also tilting.
\end{ex}

The Weyl group $W(GL(k),
T(k))$ is isomorphic to the group $S_k$ consisting of all permutations on
$k$ symbols.

More generally, one can
describe some fragment of the representation theory of any group
$GL({\bf d})$.
A maximal torus of the group $GL({\bf d})$ is $T({\bf d}
)=T(d_1)\times\ldots\times T(d_n)$. Respectively, $B(
{\bf d})=B(d_1)\times\ldots\times B(d_n)$ is a Borel subgroup and
then $B^{-}({\bf d}
)=B^{-}(d_1)\times\ldots\times B^{-}(d_n)$.
The characters of the group $T({\bf d})$ are
collections $\bar{\lambda}
=(\lambda_1,\ldots ,\lambda_n)$, where each $\lambda_i$ is a character
of the corresponding torus $T(d_i)$, $i=1,2,\ldots , n$. It is obvious that the
root data of $GL({\bf d})$ is the
direct product of the root data of the groups $GL(d_i)$. In
particular, $X(T({\bf d}))^+$ coincides with $X(T(d_1))^+\times\ldots\times
X(T(d_n))^+$.
Moreover, for any weight $\bar{\lambda}\in X(T({\bf d}))^+$ we have an
isomorphism $\bigtriangledown_{{\bf d}}(\bar{\lambda})\cong
\bigtriangledown(\lambda_1) \otimes\ldots\otimes \bigtriangledown(\lambda_n)$
and $\bigtriangleup_{{\bf d}}(\bar{\lambda})\cong
\bigtriangleup(\lambda_1)\otimes \ldots\otimes \bigtriangleup(\lambda_n)$.
Therefore, if all $\lambda_i$ are ordered partitions we see that
$\bigtriangledown_{{\bf d}}(\bar{\lambda})\cong L_{\tilde{\lambda}_1}(E_1)
\otimes\ldots L_{\tilde{\lambda}_n}(E_n)$.
The Weyl group $W(GL({\bf d}), T({\bf d}))$ is the
direct product of the Weyl groups of all factors $GL(E_i)$. In particular,
we have $\bar{\lambda}^*=(\lambda_1^*,\ldots , \lambda_n^*)$.

Consider dimensional vectors ${\bf t}(1)$, ${\bf t}(2)$ such that
${\bf d}(1)
\geq {\bf d}(2)$. We define the Schur functor $d_{{\bf d}(1),
{\bf d}(2)}$ by the following rule. For any $GL({\bf d}(1))$-module $V$ we
put $d_{{\bf d}(1), {\bf d}(2)}(V)=\sum_{\bar{\mu}\in L}
V_{\bar{\mu}}$.
The set $L$ consists of all $\bar{\mu}=(\mu_1,\ldots
,\mu_n)$ such that for any $i$ all coordinates of $\mu_i$ beginning with
$d(2)_i
+1$-th coordinate are equal to zero and $\sum_{\bar{\mu}\in X(T({\bf d}(1)))}
V_{\bar{\mu}}$
is the weight decomposition of $V$.
Identifying the group $GL({\bf d}(2))$ with
a subgroup of $GL({\bf d}(1))$ (see the subsection 1.4) we obtain that
$d_{{\bf d}(1),
{\bf d}(2)}(V)$
is a $GL({\bf d}(2))$-module. Besides, one can define a linear endomorphism
of $V$ which takes any $v=\sum_{\bar{\mu}\in X(T({\bf d}(1)))}v_{\bar{\mu}}
\in V$ to $\sum_{\bar{\mu}\in L} v_{\bar{\mu}}$. Denote this endomorphism by
the same symbol $d_{{\bf d}(1), {\bf d}(2)}$.
It is not hard to prove that if all coordinates $\lambda_i$ of $\bar{\lambda}$
are some ordered partitions then
$d_{{\bf d}(1), {\bf d}(2)}(\bigtriangleup_
{{\bf d}(1)}(\bar{\lambda}))\neq 0$ iff each "component" $\lambda_i$ has all
coordinates with numbers $\geq d(2)_i+1$ equal to zero. 
In the last case we have
$d_{{\bf d}(1), {\bf d}(2)}(\bigtriangleup_{{\bf d}(1)}(\bar{\lambda}))=
\bigtriangleup_{{\bf d}(2)}(\bar{\lambda})$. The same is valid for the
induced modules $\bigtriangledown_{{\bf d}(1)}(\bar{\lambda})$ as well as
for its simple socle. The reader can find the detailed proof in \cite{green} for
the case
$n=1$. The general case is a trivial consequence of the case $n=1$.

\subsection{ABW-filtrations}

For any vector $\lambda=
(\lambda_1,\ldots ,\lambda_s)$ with integral coordinates denote by
$\mid\lambda\mid$ its degree $\lambda_1+\ldots +\lambda_s$.
If all coordinates of $\lambda$ are non-negative integers
we denote by $\Lambda^{\lambda}(V)$ the tensor product
$\Lambda^{\lambda_1}(V)\otimes\ldots\Lambda^{\lambda_s}(V)$.

Recall the standard embedding of an exterior power $\Lambda^p(V)$ into
$V^{\otimes p}$. This map is defined by the rule

$$i_p : v_1\bigwedge\ldots\bigwedge v_p\mapsto\sum_{\sigma\in S_p}(-1)^
{\sigma} v_{\sigma(1)}\otimes\ldots\otimes v_{\sigma(p)}, v_1,\ldots, v_p\in V.
$$

\noindent Obviously, it is an $GL(V)$-equivariant. One can define
more general embedding $i_{\lambda} : \Lambda^{\lambda}(V)\rightarrow
V^{\otimes p}$, where $\lambda=(\lambda_1,\ldots ,\lambda_l)$ is any
(non-ordered) partition, $p=\mid\lambda\mid$ and $i_{\lambda}=\otimes_{1\leq q
\leq l} i_{\lambda_q}$. Denote by $p_{\lambda}$ the canonical epimorphism 
from $V^{\otimes p}$ onto
$\Lambda^{\lambda}(V)$.

Let $S^r(V\otimes W)$ be a homogeneous component of degree $r$ of the symmetric
algebra $S(V\otimes W)$, where $V, W$ are any vector spaces.
For any ordered partition $\lambda$ of degree $r$ we define the map

$$
d_{\lambda}: \Lambda^{\lambda}(V)\otimes \Lambda^{\lambda}(W)\rightarrow
S^r(V\otimes W)
$$

\noindent by $d_{\lambda} = d_{\lambda_1}\bar{\otimes} ...
\bar{\otimes} d_{\lambda_s}$, where $d_{\lambda_i}:
\Lambda^{\lambda_i}(V)\otimes \Lambda^{\lambda_i}(W)\rightarrow
S^{\lambda_i}(V\otimes W), i=1,\ldots , s,$ and
the symbol $\bar{\otimes}$ means the product map
$S^{\lambda_1}(V\otimes W)\otimes ...\otimes S^{\lambda_s}(V\otimes W)
\rightarrow S^r(V\otimes W)$.
Here, for any non-negative integer $t$ the map
$d_t :\Lambda^t(V)\otimes\Lambda^t(W)\rightarrow
S^t(V\otimes W)$ is defined by the rule

$$
d_t((v_1\wedge ...\wedge v_t)\otimes (w_1\wedge ...\wedge w_t )) =
\sum_{\sigma\in S_t}(-1)^{\sigma}v_1\otimes w_{\sigma(1)} ...
v_t\otimes w_{\sigma(t)},
$$

$$
v_i\in V, w_i\in W, 1\leq i\leq t.
$$

Let $M_{\lambda} = \sum_{\gamma\succeq\lambda}\mbox{Im} d_{\gamma}$
and $\dot{M}_{\lambda} = \sum_{\gamma\succ\lambda}\mbox{Im} d_{\gamma}$.
The symbol $\succeq$ means the lexicographical order from left to 
right on
the set of partitions. The $GL(V)\times GL(W)$-module $S^r(V\otimes W)$ has the
filtration

$$
0\subseteq M_{(r)}\subseteq M_{(r - 1,1)}\subseteq ...\subseteq
M_{\underbrace{(1,....,1)}_r} = S^r(V\otimes W)
$$

\noindent with quotients

$$
M_{\lambda}/\dot{M}_{\lambda}\cong L_{\lambda}(V)\otimes L_{\lambda}(W),
$$

\noindent where $L_{\lambda}(V)$ is the {\it Schur} module (see \cite{ak}).
We call this filtration Akin-Buchsbaum-Weyman filtration or briefly,
an ABW-filtration.

\begin{rem}
All these statements remain the same if we replace the field $K$ by any
commutative ring $R$ and require that both $V, W$ are free $R$-modules.
The functor $V\rightarrow L_{\lambda}(V)$ is universally free, that is, 
$L_{\lambda}(V)$ is a free $R$-module and commutes with change of the base
ring $R$ \cite{ak}. In particular, for any homomorphism $R\rightarrow
R'$ the functor $R'\otimes_R -$ takes ABW-filtrations to ABW-filtrations.    
\end{rem}

We define a {\it superpartition},
say
$\bar{\lambda}=(\lambda_1,
\ldots , \lambda_n)$, where $\lambda_i=(\lambda_{i1},\ldots ,
\lambda_{i, s_i}), \\
\lambda_{i1}\geq\ldots\geq\lambda_{i, s_i}\geq 0, i=1,\ldots , n$.
By definition, $\mid\bar{\lambda}\mid=\mid\lambda_1\mid+\ldots +\mid\lambda_n
\mid$.
One can endow the space

$$\Lambda^{\bar{\lambda}}({\bf f})=\prod_{i=1}^{i=n}\otimes (\Lambda^{\lambda_i}
(V_i))=\prod_{i=1}^{i=n}\otimes (\prod_{j=1}^{j=s_i}\otimes\Lambda^{\lambda_{ij}}
(V_i))$$

\noindent with a $GL({\bf f})$-module structure, where ${\bf f}=
(f_1,\ldots , f_n)$, $\dim V_i=f_i, 1\leq i\leq n$.
To be precise, each factor $GL(V_i)$
of the group $GL({\bf f})$ acts on the corresponding tensor product
$\prod_{j=1}^{j=s_i}\otimes\Lambda^{\lambda_{ij}}(V_i)$ diagonally.

It is not hard to prove that $\tilde{\bar{\lambda}}=(\tilde{\lambda}_1,\ldots ,
\tilde{\lambda}_n)$ is a highest weight of the $GL({\bf f})$-module $\Lambda^
{\bar{\lambda}}({\bf f})$. Moreover, its multiplicity is equal to $1$.
Since $\Lambda^{\bar{\lambda}}({\bf f})$ is a tilting $GL({\bf f})$-module
there are good and Weyl filtrations of this module such that the last quotient
of the first filtration (respectively the first non-zero member of the second 
one) is isomorphic to $\bigtriangledown_{{\bf f}}(\tilde{\bar{\lambda}})$
(respectively to $\bigtriangleup_{{\bf f}}(\tilde{\bar{\lambda}})$)
\cite {don4, don6, zub3, zub4}.

Denote by $R_{{\bf f}}(\bar{\lambda})$ the kernel of the corresponding
epimorphism

$$\Lambda^{\bar{\lambda}}({\bf f})\rightarrow
\bigtriangledown_{{\bf f}}(\tilde{\bar{\lambda}})$$

\noindent and by $S_{{\bf f}}(\bar{\lambda})$ the cokernel of the inclusion

$$\bigtriangleup_{{\bf f}}(\tilde{\bar{\lambda}})\rightarrow
\Lambda^{\bar{\lambda}}({\bf f}).$$

The $GL({\bf f})$-modules $R_{{\bf f}}(\bar{\lambda})$ and
$S_{{\bf f}}(\bar{\lambda})$ have GF and WF respectively. Moreover,
the module $R_{{\bf f}}(\bar{\lambda})$ and the inclusion of
$\bigtriangleup_{{\bf f}}(\tilde{\bar{\lambda}})$ are uniquely defined
(\cite{zub4}, Proposition 1.1). We have the short exact
sequence

$$0 \rightarrow S_{{\bf f}}(\bar{\lambda})^* \rightarrow
\Lambda^{\bar{\lambda}}({\bf f})^*\rightarrow \bigtriangleup_{{\bf f}}(\tilde
{\bar{\lambda}})^*\rightarrow 0.$$

By definition, $\bigtriangleup_{{\bf f}}(\tilde{\bar{\lambda}})^*\cong
\bigtriangledown_{{\bf f}}(\tilde{\bar{\lambda}}^*)$.
The unique highest weight of the module $\Lambda^{\bar{\lambda}}({\bf f})^*
\cong \Lambda^{\lambda_1}(V_1^*)\otimes\ldots\otimes\Lambda^{\lambda_n}(V_n^*)$
is equal to $\tilde{\bar{\lambda}}^*$ and since $\Lambda^{\bar{\lambda}}({\bf
f})$ is a tilting module
we get that $S_{{\bf f}}(\bar{\lambda})^*$ is uniquely defined by
the same Proposition 1.1 from \cite{zub4}.

Let us consider another group $GL({\bf g}), {\bf g}=(g_1,\ldots , g_m)$ and
some superpartition $\bar{\mu}=(\mu_1,\ldots , \mu_m),
i=1, \ldots , m$. We have the short exact
sequence of $GL({\bf f})\times GL({\bf g})$-modules

$$0\rightarrow D_{{\bf f}, {\bf g}}(\bar{\lambda}, \bar{\mu})\rightarrow
\Lambda^{\bar{\lambda}}({\bf f})\otimes\Lambda^{\bar{\mu}}({\bf g})^*
\rightarrow \bigtriangledown_{{\bf f}}(\tilde{\bar{\lambda}})\otimes
\bigtriangledown_{{\bf g}}(\tilde{\bar{\mu}}^*)
\cong \bigtriangledown_{{\bf f}}(\tilde{\bar{\lambda}})\otimes
\bigtriangleup_{{\bf g}}(\tilde{\bar{\mu}})^*
\rightarrow 0,$$

\noindent where

$$D_{{\bf f}, {\bf g}}(\bar{\lambda}, \bar{\mu})=
R_{{\bf f}}(\bar{\lambda})\otimes\Lambda^{\bar{\mu}}({\bf g})^* +
\Lambda^{\bar{\lambda}}({\bf f})\otimes S_{{\bf g}}(\bar{\mu})^* .$$

\begin{pr}
The kernel
$D_{{\bf f}, {\bf g}}(\bar{\lambda}, \bar{\mu})$ is uniquely defined and has
GF.
\end{pr}

Proof. Use the same Proposition 1.1 from \cite{zub4} and
Proposition 2.3 from \cite{dz}.

\begin{rem}
Notice that $L_{\lambda}(V^*)\cong\bigtriangleup(\tilde{\lambda})^*$
\cite{zub1}
as a $GL(V)$-module. In particular, the $GL({\bf g})$-module
$\bigtriangleup_{{\bf g}}(\tilde{\bar{\mu}})^*$ is isomorphic to
$L_{\tilde{\lambda}_1}(U_1^*)\otimes\ldots L_{\tilde{\lambda}_m}(U_m^*)$,
where $\dim U_j=g_j, 1\leq j\leq m$.
\end{rem}

\begin{rem}
In notations of Remark 1.1 there is a nondegenerate pairing of 
$GL(V)$-modules $\Lambda^t(V^*)\times
\Lambda^t(V)\rightarrow R$ defined by the rule

$$<f_1\wedge\ldots\wedge f_t, v_1\wedge\ldots\wedge v_t>=
\det(f_i(v_j)), v_i\in V, f_j\in V^*, 1\leq i, j\leq t.$$

\noindent Here $V^*=\Hom_R(V, R)$ and $V$ is a free $R$-module.
Tensoring we obtain a nondegenerate pairing $\Lambda^{\lambda}(V^*)\times
\Lambda^{\lambda}(V)\rightarrow R$. Thus $\Lambda^{\lambda}(V^*)$ is 
isomorphic to $\Lambda^{\lambda}(V)^*$ as a $GL(V)$-module. 
\end{rem}

\subsection{Specializations}

For given ${\bf d}(1)\geq {\bf d}(2)$ define
an epimorphism

$$p_{{\bf t}(1), {\bf t}(2)}: K[R(Q, {\bf t}(1))]\rightarrow K[R(Q,
{\bf t}(2))]$$

\noindent by the following rule.
Take any arrow $a\in A$. Let $i(a)=i$ and $t(a)=j$. For the sake of simplicity
denote $d_{i}(s)$ and $d_{j}(s)$ by $m_s$ and $l_s$ respectively, $s=1,2$.
We know that $m_1\geq m_2$ and $l_1\geq l_2$. Then
our epimorphism maps $y_{sr}(a)$ to zero iff either $s > l_2$ or $r >
m_2$. On the remaining variables our epimorphism is the identical map.

On the other hand, one can define the isomorphism $i_{{\bf t}(2), {\bf t}(1)}$
of the variety $R(Q, {\bf t}(2))$ onto a closed subvariety of
$R(Q, {\bf t}(1))$ by the dual rule, that is the epimorphism defined above is
the comorphism $i_{{\bf t}(2), {\bf t}(1)}^*$.

By almost the same way as $i_{{\bf t}(2), {\bf t}(1)}$ one can define the
isomorphism $j_{{\bf t}(2), {\bf t}(1)}$ of the group $H({\bf t}(2))$ onto a
closed subgroup of the group
$H({\bf t}(1))$ just bordering any invertible $d_i(2)\times d_i(2)$
matrix by the $d_i(1)-d_i(2)$ additional rows and columns which are zero
outside of the diagonal tail of length $d_i(1)-d_i(2)$. The entries on this
diagonal tail should be 1's.
It is not hard to check that
$i_{{\bf t}(2), {\bf t}(1)}(\phi^g)=
i_{{\bf t}(2), {\bf t}(1)}(\phi)^{j_{{\bf t}(2), {\bf t}(1)}(g)}$
for any $g\in H({\bf t}(2))$ and
$\phi\in R(Q, {\bf t}(2))$. The analogous equation is valid for the epimorphism
$p_{{\bf t}(1), {\bf t}(2)}$.

\subsection{Young subgroups}

Decompose an interval $[1, k]=\{1,\ldots , k\}$ into some disjoint
subsets, say
$[1, k]=\bigsqcup_{1\leq j\leq m}T_j$. Define the {\it Young subgroup}
$S_{{\cal T}}=S_{T_1}\times\ldots\times S_{T_m}$ of the group $S_k$
as  the subgroup which consists of all
permutations $\sigma\in S_k$ such that $\sigma(T_j)=T_j, 1\leq j\leq m$.
By definition, $S_T=\{\sigma\in S_k\mid\sigma(T)=T, \forall j\not\in T \
\sigma(j)=j\}$ for arbitrary subset $T$. The subsets $T_1,\ldots , T_m$ are 
said to be the {\it layers} of the group $S_{{\cal T}}$ \cite{zub1, zub4}.

The group $S_{{\cal T}}$ can be introduced in
other way. In fact, let $f$ be a map from $[1,k]$ onto $[1, m]$ defined by
the rule $f(T_j)=j, j=1,\ldots , m$. Then $S_{{\cal T}}=
\{\sigma\in S_k\mid f\circ\sigma=f\}$. Sometimes we will denote 
$S_{{\cal T}}$ by $S_f$.

For a given superpartition
$\bar{\lambda}=(\lambda_1,
\ldots , \lambda_n)$, where $\lambda_i=(\lambda_{i1},\ldots ,
\lambda_{i, s_i}),
\lambda_{i1}\geq\ldots\geq\lambda_{i, s_i}\geq 0, i=1,\ldots , n$, denote by
$S_{\bar{\lambda}}$ the Young subgroup of $S_{\mid\bar{\lambda}\mid}$
corresponding to the decomposition of $[1, \mid\bar{\lambda}\mid]$ into
sequential subintervals of lengths $\lambda_{11},\ldots , \lambda_{1, s_1},
\lambda_{21},\ldots ,\lambda_{2, s_2}, \ldots$.

For any group $G$ and its subgroup $H$ we denote by $G/H$ some fixed
representative set of the left cosets of $H$ if it does not lead to confusion.
For any $g\in G$ denote by $\bar{g}\in G/H$ the representative of the
left coset $gH$.

\subsection{Schur duality and a lemma}

Let $V$ be a projective module over a commutative ring $R$ such that, if
$f(x)\in R[x]$ is a polynomial which vanishes on $R$ then $f(x)$ is 
identically zero. 
We have a ring homomorphism 
$\psi : R[S_d]\rightarrow\End_{GL(V)}(V^{\otimes d})$ defined by the rule:
$\sigma\mapsto\tilde{\sigma}, \sigma\in S_d$, where $\tilde{\sigma}(v_1
\otimes\ldots\otimes v_d)=v_{\sigma^{-1}(1)}\otimes\ldots\otimes
v_{\sigma^{-1}(d)}$. For the sake of simplicity we will omit the upper tilde.

\begin{ther}(\cite{pr})
The homomorphism $\psi$ is surjective.
If $V$ is a free module of rank $p$ then 
the kernel $I_{p+1}$ of this epimorphism is not equal to
zero iff $d > p$ and in this case it is generated (as a
two-sided ideal) by the element $\sum_{\tau\in S_{p+1}}(-1)^{\tau}\tau$, where
$S_{p+1}=S_{[1, p+1]}$.
\end{ther}

Let $R$ be a principal ideal domain of odd or zero characteristic and 
$V, W$ are free
$R$-modules of finite ranks. For given partitions $\lambda, \mu$ of 
degree $r$ one can define an inclusion $\Phi_{\lambda, \mu}$ of 
$\Hom_R(\Lambda^{\lambda}(V), \Lambda^{\mu}(W))$ into $\Hom_R(V^{\otimes r},
W^{\otimes r})$ by the rule $\phi\mapsto i_{\mu}\phi p_{\lambda}, \phi\in
\Hom_R(\Lambda^{\lambda}(V), \Lambda^{\mu}(W))$. 

\begin{lm}
The image of $\Hom_R(\Lambda^{\lambda}(V), \Lambda^{\mu}(W))$ in
$\Hom_R(V^{\otimes r}, W^{\otimes r})$ coincides with 
$\{\phi\in\Hom_R(V^{\otimes r}, W^{\otimes r})\mid \forall\tau_1\in 
S_{\lambda}, \forall\tau_2\in S_{\mu}, \tau_2\phi\tau_1=(-1)^{\tau_1}
(-1)^{\tau_2}\phi\}$.  
\end{lm}

Proof. Denote the submodule $\{\phi\in\Hom_R(V^{\otimes r}, W^{\otimes r})
\mid \forall\tau_1\in S_{\lambda}, \forall\tau_2\in S_{\mu}, \tau_2\phi\tau_1
=(-1)^{\tau_1} (-1)^{\tau_2}\phi\}$ by $M$. By Remark 1.3 we have the 
standard isomorphisms

$$\Hom_R(\Lambda^{\lambda}(V), \Lambda^{\mu}(W))\cong\Lambda^{\lambda}
(V^*)\otimes\Lambda^{\mu}(W), \Hom_R(V^{\otimes r}, W^{\otimes r})\cong
(V^*)^{\otimes r}\otimes W^{\otimes r}.$$     

\noindent 
Then $\Phi_{\lambda, \mu}$ can be identified with $i_{\lambda}\otimes 
i_{\mu}$. The groups $S_{\mu}$ and $S_{\lambda}$ act on $(V^*)^{\otimes r}
\otimes W^{\otimes r}$ in obvious way. It remains to notice that for any 
free $R$-module $U$ and partition $\chi$ (of $r$) we have
$i_{\chi}(\Lambda^{\chi}(U))=\{x\in U^{\otimes r}\mid\forall\tau\in S_{\chi},
\tau(x)=(-1)^{\tau} x\}$. In fact, $\mbox{Im} \Phi_{\lambda, \mu}=   
\mbox{Im} i_{\lambda}\otimes i_{\mu}=M$.

\section{Auxiliary computations}

\subsection{Generators, free invariant algebras and relations}

We start with some simplification of the space $R(Q, {\bf t})$.
To be precise, let $a\in A$ and $i(a)=i, t(a)=j$. We have the following
possibilities:

\begin{enumerate}
\item If $W_i=E_i,
W_j=E_j$ then $H=H({\bf t})$ acts on the component $K[\Hom_K(E_i, E_j)]=
K[Y(a)]=K[y_{lt}(a)\mid 1\leq l\leq d_j, 1\leq t\leq d_i]$ by the rule
$Y(a)\mapsto g^{-1}Y(a)h, g\in GL(d_j), h\in GL(d_i)$. It can easily be
checked that
$K[Y(a)]\cong S(E_j^*\otimes E_i)$ and this isomorphism of $GL(d_j)\times
GL(d_i)$-modules is defined by the rule $y_{lt}(a)\longleftrightarrow
e_l^*\otimes f_t$, where $e_1,\ldots , e_{d_j}$ and $f_1,\ldots , f_{d_i}$
are some fixed bases of the spaces $E_j$ and $E_i$ respectively.
The basis $e_1^*,\ldots , e_{d_j}^*$ is the dual relative to $e_1,\ldots ,
e_{d_j}$.
\item If $W_i=E_i, W_j=E_j^*$ then $K[Y(a)]\cong S(E_j\otimes E_i)$
with respect to the identification $y_{lt}(a)\longleftrightarrow
e_l\otimes f_t$. In other words, $H$ acts on $Y(a)$ by 
$Y(a)\mapsto g^tY(a)h$.

Other cases are listed without any comments.

\item $W_i=E_i^*, W_j=E_j$, $K[Y(a)]\cong S(E_j^*\otimes E_i^*),
y_{lt}(a)\longleftrightarrow e_l^*\otimes f_t^*$,
$Y(a)\mapsto g^{-1}Y(a)(h^t)^{-1}$.

\item $W_i=E_i^*, W_j=E_j^*$, $K[Y(a)]\cong S(E_j\otimes E_i^*),
y_{lt}(a)\longleftrightarrow e_l\otimes f_t^*$,
$Y(a)\mapsto g^tY(a)(h^t)^{-1}$.
\end{enumerate}

\begin{lm} Up to some changing of $Q$ one can eliminate 
the fourth case. 
\end{lm}

Proof. By the definition there should be some $q, q'\in\Omega$ such that
$i=j_q, j=j_{q'}$. Redefine the maps $i, t$ on any arrow $a$
which goes from $i$ to $j$
by : $i'(a)=i_{q'},
t'(a)=i_q$ and $i', t'$ coincide with $i, t$ on the remaining arrows.
We get a new quiver $Q'$.

Let us consider the representation space
of this new quiver of the same dimension ${\bf t}$. It is clear that this space
can be produced from $R(Q, {\bf t})$ by
replacing all summands $\Hom_K(E_i^*, E_j^*)$
by $\Hom_K(E_{j}, E_{i})$.

The group $H$ remains the same.
Moreover, the algebra $K[R(Q, {\bf t})]$ is isomorphic to
$K[R(Q', {\bf t})]$. To be precise, we map each $y_{lt}(a)$
to $z_{tl}(a)$, where $Z(a)=Z_{{\bf t}}(a), i'(a)=i_{q'}, t'(a)=i_q$.
The remaining generators of
$K[R(Q, {\bf t})]$ and $K[R(Q', {\bf t})]$ coincide with one another.
It can easily be checked that this isomorphism is $H$-equivariant.
After repeating this procedure as many times as we need one can see that
the fourth case does not happen at all. The lemma is proved.

Decompose the arrow set $A$ into three subsets $A_i, i=1, 2, 3$, where
$A_1=\{a\in A\mid W_{i(a)}=E_{i(a)}, W_{t(a)}=E_{t(a)}\}$,
$A_2=\{a\in A\mid W_{i(a)}=E_{i(a)}, W_{t(a)}=E_{t(a)}^*\}$ and
$A_3=\{a\in A\mid W_{i(a)}=E_{i(a)}^*, W_{t(a)}=E_{t(a)}\}$.

The algebra $K[R(Q, {\bf t})]$ is isomorphic to the tensor product

$$
\prod_{1\leq k\leq 3}\otimes (\otimes_{a\in A_k}K[Y(a)])
$$

\noindent or to

$$
\prod_{1\leq k\leq 3}\otimes (\otimes_{a\in A_k}
(\oplus_{r_a} K[Y(a)](r_a)))\cong
(\prod_{a\in A_1}\otimes (\oplus_{r_a}S^{r_a}(E^*_{t(a)}\otimes E_{i(a)}))
\otimes$$

$$(\prod_{a\in A_2}\otimes (\oplus_{r_a}S^{r_a}(E_{t(a)}\otimes
E_{i(a)}))\otimes
(\prod_{a\in A_3}\otimes (\oplus_{r_a}S^{r_a}(E^*_{t(a)}\otimes E^*_{i(a)})),
$$

\noindent as a $H$-module.

Fix a {\it multidegree} $\bar{r}=(r_a)_{a\in A}$. Sometimes we will
rewrite it as $(\bar{r}_1,\bar{r}_2,\bar{r}_3)$, where $\bar{r}_i=(r_a)_{a\in
A_i}, i=1,2,3$. Denote $\sum_{a\in A}r_a$ by $r$ and $\sum_{a\in A_i}r_a$ by
$r_i, i=1,2,3$. The $\bar{r}$-homogeneous component of the algebra
$K[R(Q,{\bf t})]$ is isomorphic to

$$(\prod_{a\in A_1}\otimes S^{r_a}(E^*_{t(a)}\otimes E_{i(a)}))\otimes
(\prod_{a\in A_2}\otimes S^{r_a}(E_{t(a)}\otimes E_{i(a)}))\otimes
(\prod_{a\in A_3}\otimes S^{r_a}(E^*_{t(a)}\otimes E^*_{i(a)})).
$$

Tensoring ABW-filtrations  of all factors in this tensor product
we see that the $\bar{r}$-homogeneous component of the algebra
$K[R(Q,{\bf t})]$ has a filtration with quotients

$$
\prod_{a\in A_1}\otimes (L_{\lambda_a}
(E_{t(a)}^*)\otimes L_{\lambda_a}(E_{i(a)}))\otimes
\prod_{b\in A_2}\otimes (L_{\mu_b}
(E_{t(a)})\otimes L_{\mu_b}(E_{i(a)}))\otimes$$
$$ 
\otimes\prod_{c\in A_3}\otimes (L_{\gamma_c}
(E_{t(c)}^*)\otimes L_{\gamma_c}(E_{i(c)}^*))
$$

\noindent as a $\prod_{1\leq l\leq 3}(\prod_{a\in A_l} (GL(d_{t(a)})\times
GL(d_{i(a)}))$-module, where by definition $\forall a\in A_1, \forall b\in
A_2, \forall c\in A_3, \mid\lambda_a \mid=r_a, \mid\mu_b\mid=r_b, \mid\gamma_c
\mid=r_c$ (see \cite{dz}, Proposition 2.3).
We enumerate the members of this filtration
by the triples (superpartitions) $\Theta=(\lambda_{A_1}, \mu_{A_2},
\gamma_{A_3})$, 
where $\lambda_{A_1}=(\lambda_a)_{a\in A_1}, \mu_{A_2}=(\mu_a)_{a\in A_2},
\gamma_{A_3}=(\gamma_a)_{a\in A_3}$, say

\begin{equation}
\ldots\subseteq M_{\Theta}({\bf t})=M_{\Theta}\subseteq\ldots .
\end{equation}

Denote by
$\Lambda_1(\Theta, {\bf t})$
and $\Lambda_2(\Theta, {\bf t})$ the spaces

$$\prod_{a\in A_1}\otimes (\Lambda^{\lambda_a}(E_{i(a)}))
\otimes \prod_{a\in A_2}\otimes (\Lambda^{\mu_a}
(E_{t(a)})\otimes \Lambda^{\mu_a}(E_{i(a)}))$$

and

$$\prod_{a\in A_1}\otimes (\Lambda^{\lambda_a}
(E_{t(a)}))\otimes
\prod_{a\in A_3}\otimes (\Lambda^{\gamma_a}
(E_{t(a)})\otimes \Lambda^{\gamma_a}(E_{i(a)}))
$$

\noindent respectively. Sometimes we will omit the indices ${\bf t}, 
{\bf d}$ or
$\Theta$ if it does not lead to confusion.

By Remark 1.3 one can identify  the dual space
$\Lambda_2(\Theta, {\bf t})^*$ with the space

$$
\prod_{a\in A_1}\otimes (\Lambda^{\lambda_a}
(E_{t(a)}^*))\otimes
\prod_{a\in A_3}\otimes (\Lambda^{\gamma_a}
(E_{t(a)}^*)\otimes \Lambda^{\gamma_a}(E_{i(a)}^*)).
$$

Arrange the tensor factors of the quotient $M_{\Theta}/\dot{M}_{\Theta}$
into groups by the following way:

$$
\prod_{a\in A_1}\otimes (L_{\lambda_a}(E_{i(a)}))
\otimes \prod_{a\in A_2}\otimes (L_{\mu_a}
(E_{t(a)})\otimes L_{\mu_a}(E_{i(a)}))\otimes$$

$$
\prod_{a\in A_1}\otimes (L_{\lambda_a}
(E_{t(a)}^*))\otimes
\prod_{a\in A_3}\otimes (L_{\gamma_a}
(E_{t(a)}^*)\otimes L_{\gamma_a}(E_{i(a)}^*)).
$$

The first factor

$$
\prod_{a\in A_1}\otimes (L_{\lambda_a}(E_{i(a)}))
\otimes \prod_{a\in A_2}\otimes (L_{\mu_a}
(E_{t(a)})\otimes L_{\mu_a}(E_{i(a)}))$$

\noindent is a $(\prod_{a\in A_1} GL(d_{i(a)}))\times (\prod_{a\in A_2}
GL(d_{t(a)})\times GL(d_{i(a)}))$-module. Denote this group by $G_1=
G_1({\bf d})$.

Analogously, the second factor

$$
\prod_{a\in A_1}\otimes (L_{\lambda_a}
(E_{t(a)}^*))\otimes
\prod_{a\in A_3}\otimes (L_{\gamma_a}
(E_{t(a)}^*)\otimes L_{\gamma_a}(E_{i(a)}^*))
$$

\noindent is a $(\prod_{a\in A_1} GL(d_{t(a)}))\times (\prod_{a\in A_3}
GL(d_{t(a)})\times GL(d_{i(a)}))$-module. Denote it by $G_2
=G_2({\bf d})$.

For any $\Theta$ we have a homomorphism of $G_1\times G_2$-modules

\begin{equation}
d_{\Theta} : \Lambda_1(\Theta)
\otimes \Lambda_2(\Theta)^*\rightarrow K[R(Q, {\bf t})](\bar{r}).
\end{equation}

The homomorphism $d_{\Theta}$
induces an epimorphism 

$$\Lambda_1(\Theta)
\otimes \Lambda_2(\Theta)^*\rightarrow M_{\Theta}/\dot{M}_{\Theta}\rightarrow 
0.$$

Denote by $\Theta_1, \Theta_2$ the superpartitions 
$(\lambda_{A_1},\mu_{A_2}, \mu_{A_2})$ and
$(\lambda_{A_1},\gamma_{A_3}, \gamma_{A_3})$
respectively.

The $G_1({\bf d})$-module

$$(\prod_{a\in A_1}\otimes L_{\lambda_a}(E_{i(a)}))
\otimes\prod_{a\in A_2}\otimes (L_{\mu_a}
(E_{t(a)})\otimes L_{\mu_a}(E_{i(a)}))$$

\noindent coincides with

$$\bigtriangledown_{((d_{i(a)})_
{a\in A_1}, (d_{t(a)})_{a\in A_2}, (d_{i(a)})_{a\in A_2})}
(\tilde{\Theta_1}).$$

Similarly, the $G_2({\bf d})$-module
$$(\prod_{a\in A_1}\otimes L_{\lambda_a}
(E_{t(a)}^*))\otimes \prod_{a\in A_3}\otimes
(L_{\gamma_a}(E_{t(a)}^*)\otimes L_{\gamma_a}(E_{i(a)}^*))$$
coincides with 
$$\bigtriangleup_{((d_{t(a)})_
{a\in A_1}, (d_{t(a)})_{a\in A_3}, (d_{i(a)})_{a\in A_3})}
(\tilde{\Theta_2})^*.$$
Slightly abusing our notations we denote these modules by 
$\bigtriangledown_{{\bf d}}(\Theta)$ and $\bigtriangleup_{{\bf d}}(\Theta)$
correspondingly.
 
Using Proposition 1.2 we obtain the uniquely defined short exact
sequence of $G_1\times G_2$-modules with GF

$$0\rightarrow D_{{\bf d}, {\bf d}}(\Theta)=D(\Theta)\rightarrow
\Lambda_1(\Theta)\otimes
\Lambda_2(\Theta)^*\rightarrow M_{\Theta}/\dot{M}_{\Theta}\rightarrow 0.$$

\noindent Here, 

$$D_{{\bf d}, {\bf d}}(\Theta)=D_{((d_{i(a)})_
{a\in A_1}, (d_{t(a)})_{a\in A_2}, (d_{i(a)})_{a\in A_2}),
((d_{t(a)})_{a\in A_1}, (d_{t(a)})_{a\in A_3}, (d_{i(a)})_{a\in A_3})}
(\Theta_1, \Theta_2).$$

To turn
to the group $H=H({\bf d})$ we have to replace the group $G_1$($G_2$)
by some subgroup. Indeed, represent, say $G_1$, as
$\times_{i\in V}GL(d_i)^{w_i}$, where $w_i$ is the number of factors of
$G_1$ coinciding with $GL(d_i), i\in V$. The next step is to replace any
subproduct $GL(d_i)^{w_i}, i\in V_{ord},$ or $GL(d_q)^{v_q}, q\in\Omega$,
where
$v_q=w_{i_q}+w_{j_q}$,
by the corresponding diagonal subgroup.

Using Theorem 1.1(3) we obtain that any $G_i$-module has GF
(respectively -- any $G_i$-module has WF) retains this property under the
restriction to the group $H$, $i=1,2$. Referring to Theorem 1.1(1)
we get

\begin{pr}
The short sequence

$$0\rightarrow D(\Theta)^{H}\rightarrow
(\Lambda_1(\Theta)\otimes
\Lambda_2(\Theta)^*)^{H}\rightarrow Z(\Theta)
\rightarrow 0$$

is exact.
Here, $Z(\Theta)=(M_{\Theta}/\dot{M}_{\Theta})^H=
M_{\Theta}^H/\dot{M}_{\Theta}^H$.
\end{pr}

The same arguments show that $\Lambda_1$ and $\Lambda_2$ are tilting 
$H$-modules and all quotients of the filtration (1)  
are $H$-modules with GF.

One can rewrite the exact sequence from Proposition 2.1 as

\begin{equation}
0\rightarrow D(\Theta)^{H}\rightarrow
\Hom_H(\Lambda_2(\Theta), \Lambda_1(\Theta))\rightarrow
Z(\Theta)\rightarrow 0 .
\end{equation}

Sometimes, if it is necessary to indicate that the original representation 
space has
dimension ${\bf t}$, we write $Z_{{\bf t}}(\Theta)$.

\begin{ther}(\cite{don2})
The epimorphism $p_{{\bf t}(1), {\bf t}(2)} : K[R(Q,{\bf t}(1))]\rightarrow
K[R(Q,{\bf t}(2))]$ induces the epimorphism $\phi_{{\bf t}(1), {\bf t}(2)} :
J(Q,{\bf t}(1))\rightarrow J(Q,{\bf t}(2))$.
\end{ther}

Proof. For given pair of compatible dimension vectors ${\bf d}(1)\geq
{\bf d}(2)$
one can define at least three Schur functors $d, d_1, d_2$ for the groups
$H, G_1, G_2$ correspondingly.
Nevertheless, it is not hard to see that the action of the Schur functor $d$
coincides with the actions of both functors $d_i, i=1,2$
on the short exact sequences

$$0\rightarrow R_{{\bf d}(1)}(\Theta)\rightarrow \Lambda_1({\bf t}(1))
\rightarrow\bigtriangledown_{{\bf d}(1)}(\Theta)\rightarrow 0$$

\noindent and

$$0\rightarrow \bigtriangleup_{{\bf d}(1)}(\Theta)\rightarrow
\Lambda_2({\bf t}(1))\rightarrow S_{{\bf d}(1)}(\Theta)\rightarrow 0.$$

\noindent It is obvious for modules $\Lambda_1(\Theta, {\bf t}(1)), 
\Lambda_2(\Theta,
{\bf t}(1)),
\bigtriangledown_{{\bf d}(1)}(\Theta),
\bigtriangleup_{{\bf d}(1)}(\Theta)$ and it follows for
$R_{{\bf d}(1)}(\Theta), S_{{\bf d}(1)}(\Theta)$ since all Schur functors are
exact \cite{green}.
Moreover, one can identify the exact sequences

$$0\rightarrow d(R_{{\bf d}(1)}(\Theta))\rightarrow d(\Lambda_1({\bf t}(1)))
\rightarrow
d(\bigtriangledown_{{\bf d}(1)}(\Theta))\rightarrow 0$$

\noindent and

$$0\rightarrow d(\bigtriangleup_{{\bf d}(1)}(\Theta))\rightarrow
d(\Lambda_2({\bf t}(1)))\rightarrow d(S_{{\bf d}(1)}(\Theta))\rightarrow 0$$

\noindent with

$$0\rightarrow R_{{\bf d}(2)}(\Theta)\rightarrow \Lambda_1({\bf t}(2))
\rightarrow
\bigtriangledown_{{\bf d}(2)}(\Theta)\rightarrow0$$

\noindent and

$$0\rightarrow \bigtriangleup_{{\bf d}(2)}(\Theta)\rightarrow
\Lambda_2({\bf t}(2))\rightarrow S_{{\bf d}(2)}(\Theta)\rightarrow 0$$

\noindent respectively since $R$ and $S$ are uniquely defined in all these 
sequences.

Let $\psi: \Lambda_1\otimes \Lambda_2^*\rightarrow d(\Lambda_1)
\otimes d(\Lambda_2)^*$ be a map given by $\psi(v\otimes\alpha)=d(v)\otimes
\alpha\mid_{d(\Lambda_2)}$,
$v\in \Lambda_1, \alpha\in \Lambda_2^*$. In other words, it is the map 
$\Hom_K(\Lambda_2, \Lambda_1)\rightarrow\Hom_K(d(\Lambda_1), 
d(\Lambda_2))$ defined by $\phi\mapsto d\circ\phi\mid_{d(\Lambda_2)}$.

If $\phi\in\Hom_H(\Lambda_2, \Lambda_1)$ then $\phi(d(\Lambda_2))\subseteq d(\Lambda_1)$ since
$\phi$ commutes with the torus action. In particular, $\psi$ is the restriction
map on $\Hom_H(\Lambda_2, \Lambda_1)$. Moreover, $\psi(D_{{\bf d}(1)})\subseteq
D_{{\bf d}(2)}$.
Indeed, it is clear for the summand $R\otimes \Lambda_2^*$. Let $v\otimes\alpha\in
\Lambda_1\otimes S^*$. The space $S^*$ is identified with a subspace of $\Lambda_2^*$ by the
rule $\alpha\mapsto\alpha\circ p$, where $p$ is the epimorphism of the
$G_2$-modules $\Lambda_2\rightarrow S\rightarrow 0$. In particular, 
$p(d(\Lambda_2))=d(S)$
and $(\alpha\circ p)\mid_{d(\Lambda_2)}=\alpha\mid_{d(S)}\circ 
p\mid_{d(\Lambda_2)}$.

Consider the filtration $0\subseteq R\otimes S^*\subseteq D$
of $H$-module $D$
with quotients $R\otimes S^*$ and $(R\otimes\bigtriangleup^*)\oplus
(\bigtriangledown\otimes S^*)$. These quotients
can be identified with $\Hom_K(S, R)$ and $\Hom_K(\bigtriangleup, R)\oplus
\Hom_K(S, \bigtriangledown )$ respectively and the map $\psi$ induces
the maps

$$\Hom_K(S, R)\rightarrow\Hom_K(d(S), d(R)),
\Hom_K(\bigtriangleup, R)\rightarrow\Hom_K(d(\bigtriangleup), d(R)),$$
$$\Hom_K(S, \bigtriangledown )\rightarrow\Hom_K(d(S), d(\bigtriangledown) ).$$

All these arguments show that we have the following
commutative diagram

$$\begin{array}{ccccccccc}
0 & \rightarrow & D & \rightarrow & \Hom_K(\Lambda_2, \Lambda_1) & \rightarrow
&\bigtriangledown\otimes \bigtriangleup^* & \rightarrow & 0 \\
& & \downarrow & & \downarrow & & \downarrow & &  \\
0 & \rightarrow & d(D) & \rightarrow & \Hom_K(d(\Lambda_2), d(\Lambda_1)) &
\rightarrow & d(\bigtriangledown)\otimes d(\bigtriangleup)^* & \rightarrow & 0.
\end{array}
$$

If we identify the last right members of the horizontal sequences with the
corresponding
quotients of the filtrations of $K[R(Q, {\bf t}(1))](\bar{r})$
and $K[R(Q, {\bf t}(2))](\bar{r})$ respectively then the last right vertical
arrow is
induced by the epimorphism $p_{{\bf t}(1), {\bf t}(2)}$.
Indeed, the map $d=d_{{\bf d}(1), {\bf d}(2)}$
takes a basis vector of $\Lambda_1=\Lambda_1({\bf t}(1))$ or $\Lambda_2=\Lambda_2({\bf t}(1))$ to zero if its
record contains at least one
vector $e_j^{(i)}$ or $(e_j^{(i)})^*$, where $j\geq d(2)_i + 1$ and
$e_1^{(i)},\ldots , e_{d(1)_i}^{(i)}$ is a fixed basis of $E_i$, $1\leq i\leq
n$. It remains to apply the rule of the identification of the algebra
$K[R(Q, {\bf t})]$ with the corresponding symmetric algebra.

Finally, we have the following commutative diagram

\begin{equation}
\begin{array}{ccccccccc}
0 & \rightarrow & D_{{\bf d}(1),{\bf d}(1)}^{H({\bf t}(1))} & \rightarrow &
\Hom_{H({\bf t}(1))}(\Lambda_2({\bf t}(1)), \Lambda_1({\bf t}(1))) & 
\rightarrow & Z_{{\bf t}(1)} & \rightarrow & 0 \\
& & \downarrow & & \downarrow & & \downarrow & &  \\
0 & \rightarrow & D_{{\bf d}(2),{\bf d}(2)}^{H({\bf t}(2))} & \rightarrow &
\Hom_{H({\bf t}(2))}(\Lambda_2({\bf t}(2)), \Lambda_1({\bf t}(2))) &
\rightarrow & Z_{{\bf t}(2)} & \rightarrow & 0.
\end{array}
\end{equation}

Repeating word by word the proof of Proposition 1 from \cite{zub1} (and
using Lemma 1.1 from \cite{zub4} as well) we see that all vertical
arrows in
the last diagram are epimorphisms. This concludes the proof.

We have an inverse spectrum of algebras:

$$
\{ J(Q, {\bf t}),\phi_{{\bf t}(1), {\bf t}(2)}\mid {\bf d}
(1)\geq {\bf d}(2)\}.
$$

Moreover, because of epimorphisms $\phi_{{\bf t}(1), {\bf t}(2)}$ are
homogeneous we have the countable set of spectrums:

$$
\{ J(Q, {\bf t})(r), \phi_{{\bf t}(1), {\bf t}(2)} \mid {\bf d}
(1)\geq {\bf d}(2)\}, r=0,1,2,\ldots\}. $$

The inverse limit of $r$-th spectrum denote by $J(Q)(r)$. It is clear that
$J(Q)=\oplus_{r\geq 0}J(Q)(r)$ can be endowed with an algebra structure in
obvious way.
The algebra $J(Q)$ is called a {\it free} invariant algebra of
mixed representations of the quiver $Q$.

\begin{rem}
From the geometrical point of view we have a commutative diagram

$$\begin{array}{ccc}
R(Q, {\bf t}(2)) & \rightarrow & R(Q, {\bf t}(2))/H({\bf t}(2))\\
\downarrow &  & \downarrow \\
R(Q, {\bf t}(1)) & \rightarrow & R(Q, {\bf t}(1))/H({\bf t}(1))
\end{array}$$

\noindent which is dual to

$$\begin{array}{ccc}
J(Q, {\bf t}(2)) & \rightarrow & K[R(Q, {\bf t}(2))]\\
\uparrow &  & \uparrow \\
J(Q, {\bf t}(1)) & \rightarrow & K[R(Q, {\bf t}(1))].
\end{array}$$

In the first diagram horizontal sequences are categorical quotients with
respect to the
corresponding reductive group actions and vertical arrows are isomorphisms
onto closed subvarieties. The algebra $J(Q)$ can be regarded as a coordinate
algebra of an infinitely dimensional  variety which is the direct limit of
varieties  $\Spec(J(Q,{\bf t}))=R(Q, {\bf t})/H({\bf t})$ or as an invariant
algebra $K[R(Q)]^{H(Q)}$, where $K[R(Q)]$ is the {\it homogeneous} inverse
limit of the algebras $K[R(Q,{\bf t})]$ defined by the same way as above and
$H(Q)$ is the direct limit of the groups $H({\bf t})$.
\end{rem}

It is clear that any $J(Q, {\bf t})$ is an epimorphic image of $J(Q)$.
Denote the kernel of this epimorphism by $T(Q,{\bf t})$.

\begin{rem}
It is not necessary to consider the algebra $J=J(Q)$
as the inverse limit over all compatible dimensional vectors
${\bf t}$. One can replace the set of all dimensional vectors by any cofinal
subset.
For example, we can take $\{{\bf N} = (T_1,\ldots , T_n)\mid N\geq 2\}$,
where $T_i=N$ iff $t_i=d_i$ otherwise $T_i=N^*$.
\end{rem}

From now on we
suppose that ${\bf t}(1)={\bf N}$ and ${\bf t}(2)={\bf t}$, where the number
$N$ is sufficiently large, say $N\geq r$. Denote by ${\bf d}$ the underlying
vector of ${\bf t}$.

Finally, denote the image $d_{\Theta}(\phi)$ of any 
$\phi\in\Hom_{H({\bf t})}(
\Lambda_2({\bf t}),
\Lambda_1({\bf t}))$ in the homogeneous component $K[R(Q, {\bf t})](\bar{r})$
by $c(\phi)$. 

\begin{lm}
If $\Hom_{H({\bf t})}(\Lambda_2(\Theta), \Lambda_1(\Theta))\neq 0$ then
$r_2=r_3=\mid\mu_{A_2}\mid=\mid\gamma_{A_3}\mid$. Moreover, the following 
conditions are satisfied :

\begin{enumerate}
\item $\forall i\in V_{ord}, \sum_{a\in A, t(a)=i}r_a =
\sum_{a\in A, i(a)=i} r_a=p_a.$
\item $\forall q\in\Omega, \sum_{a\in A, t(a)=i_q} r_a +\sum_{a\in A, i(a)=j_q}
r_a=\sum_{a\in A, i(a)=i_q} r_a +\sum_{a\in A, t(a)=j_q} r_a=p_q.$
\end{enumerate} 
\end{lm}

Proof.
It is clear that the group $H({\bf N})$ contains the diagonal subgroup which
is isomorphic to $GL(N)$. The requirement
$\Hom_{H({\bf t})}(\Lambda_2(\Theta), \Lambda_1(\Theta))\neq 0$ 
implies $\Hom_{H({\bf N})}
(\Lambda_2({\bf N}), \Lambda_1({\bf N}))\neq 0$ for sufficiently large $N$.
Therefore, we obtain that $\Hom_{GL(N)}(\Lambda_2({\bf N}),
\Lambda_1({\bf N}))\neq 0$ and  
the degrees of the polynomial $GL(N)$-modules (see \cite{green} for
definitions)
$\Lambda_1(\Theta)$ and
$\Lambda_2(\Theta)$ should be the same. In other words,
$r_2=r_3=\mid\mu_{A_2}\mid=
\mid\gamma_{A_3}\mid$.

The space $\Hom_{H({\bf t})}(\Lambda_2(\Theta), \Lambda_1(\Theta)$ can be 
represented as

$$\otimes_{i\in V_{ord}}\Hom_{GL(d_i)}(\otimes_{a\in A, t(a)=i}
\Lambda^{\chi_a}(E_i), \otimes_{a\in A, i(a)=i}\Lambda^{\chi_a}(E_i))
$$

$$\otimes_{q\in\Omega}\Hom_{GL(d_q)}((\otimes_{a\in A, t(a)=i_q}
\Lambda^{\chi_a}(E_{i_q}))
\otimes (\otimes_{a\in A, i(a)=j_q}
\Lambda^{\chi_a}(E_{j_q})),$$

$$(\otimes_{a\in A, i(a)=i_q}\Lambda^{\chi_a}(E_{i_q}))
\otimes (\otimes_{a\in A, t(a)=j_q}\Lambda^{\chi_a}(E_{j_q}))),$$

\noindent where $\chi_a$ is equal to $\lambda_a, \mu_a$ or $\gamma_a$ if
$a\in A_1, a\in A_2$ or $a\in A_3$ respectively.
All tensor multipliers are not equal to zero. It remains to compare the 
degrees of all modules in the corresponding groups of homomorphisms. 

Denote $\mid\lambda_{A_1}\mid$ by $t$. Then $r=t+2s$.

As in \cite{zub4} we extend the set of matrix variables $\{Y(a)\mid a\in A\}$
in the following way. Replace each $Y(a)$ by some new set of matrices
having the same size as $Y(a)$. The cardinality of this set is equal to $r_a$.
Simultaneously, we replace each arrow $a$ by $r_a$ new arrows with the same
origin
and end as $a$ and set them in one-to-one correspondence with these new 
matrices.
So we get a new quiver $\hat{Q}$.
The vertex set of $\hat{Q}$ coincides with $V$ but the arrow set $\hat{A}$
can be different from $A$.

Take any linear order on $A$. Denote this order by usual
symbol $<$. We enumerate arrows of the quiver
$\hat{Q}$
by numbers $1,\ldots , r$. One can assume that for any $a\in A$ the
corresponding set
of new arrows is enumerated by the numbers from the segment $[\dot{a}, a]=
[\sum_{b<a}r_b+1,\sum_{b\leq a}r_b]$.
We obtain some {\it specialization} $f: [1, r]=\hat{A}
\rightarrow A$ defined by $f(j)=a$ iff $j\in [\dot{a}, a], a\in A$.

In the same way one can define the specialization $Y(j)\mapsto Y(a)$ iff
$j\in [\dot{a}, a], a\in A$. Denote the last specialization by the same symbol
$f$.

Without loss of generality it can be assumed that $\forall a\in\hat{A}_1,
b\in\hat{A}_2, c\in\hat{A}_3, a < b < c$. Thus 
$\hat{A}_1=[1,t], \hat{A}_2=[t+1,t+s], \hat{A}_3=[t+s+1, r]$.
Moreover,
$f([1, t])=A_1$,
$f([t+1, s+t])=A_2$ and $f([s+t+1, r])=A_3$. It is clear that
$i(j)=i$ or $t(j)=i$ iff $i(f(j))=i$ or $t(f(j))=i$ respectively,
$j\in\hat{A}=[1,\ldots , r], i\in V$.
Set

$$T(i)=
\{j\in\hat{A}\mid t(j)=i\}, I(i)=\{j\in\hat{A}\mid i(j)=i\}, i\in V_{ord}.$$

Analogously, write

$$T(q)=\{j\in\hat{A}\mid t(j)=i_q \ \mbox{or} \ i(j)=j_q\},
I(q)=\{j\in\hat{A}\mid i(j)=i_q \ \mbox{or} \ t(j)=j_q\}, q\in\Omega.$$

It is obvious that
$p_i=\mid T(i)\mid=\mid I(i)\mid$ for each $i\in V_{ord}$ and $p_q=\mid
T(q)\mid=\mid I(q)\mid$ for all $q\in\Omega$.

We have the inclusion $\Phi_{\Theta}=\Phi_{\Theta_2, \Theta_1}$ of the space
$\Hom_{H({\bf t})}(\Lambda_2(\Theta), \Lambda_1(\Theta))$ 
into

$$\Hom_{H({\bf t})}((\otimes_{a\in A_1}
E_{t(a)}^{\otimes r_a})\otimes
(\otimes_{a\in A_3} E_{t(a)}^{\otimes r_a})\otimes (\otimes_{a\in A_3}
E_{i(a)}^{\otimes r_a}),
$$

$$(\otimes_{a\in A_1}(E_{i(a)}^{\otimes r_a})\otimes
(\otimes_{a\in A_2} E_{t(a)}^{\otimes r_a})\otimes (\otimes_{a\in A_2}
E_{i(a)}^{\otimes r_a})).
$$

Denote the last
space by $\Hom({\bf t})$.

The multilinear component of degree $r$ of the algebra
$J(\hat{Q},{\bf t})$ is isomorphic to

$$\Hom_{H({\bf t})}((\otimes_{a\in \hat{A}_1}
E_{t(a)})\otimes
(\otimes_{a\in \hat{A}_3}E_{t(a)})\otimes (\otimes_{a\in \hat{A}_3}E_{i(a)}),
$$

$$(\otimes_{a\in \hat{A}_1}(E_{i(a)})\otimes
(\otimes_{a\in \hat{A}_2} E_{t(a)})\otimes (\otimes_{a\in \hat{A}_2}
E_{i(a)})).
$$

It is clear that this space coincides with $\Hom({\bf t})$.

We identify the space
$\Hom_{H({\bf N})}(\Lambda_2({\bf N}),\Lambda_1({\bf N}))$ with its image in
$\Hom({\bf N})$.

\begin{lm}
If both modules
$\Lambda_2({\bf t}),\Lambda_1({\bf t})$ are not equal to zero then 
the kernel of the map

$$\Hom_{H({\bf N})}(\Lambda_2({\bf N}),\Lambda_1({\bf N}))\rightarrow
\Hom_{H({\bf t})}(\Lambda_2({\bf t}),\Lambda_1({\bf t}))$$

\noindent is the intersection of
$\Hom_{H({\bf N})}(\Lambda_2({\bf N}),\Lambda_1({\bf N}))$
with the kernel of the epimorphism
$\Hom({\bf N})\rightarrow \Hom({\bf t})$.
\end{lm}

Proof. It is sufficient to look at the following commutative diagram

$$
\begin{array}{ccccc}
0 & \rightarrow & \Hom_{H({\bf N})}(\Lambda_2({\bf N}),\Lambda_1({\bf N})) &
\rightarrow
& \Hom({\bf N})\\
 & & \downarrow & & \downarrow \\
0 & \rightarrow & \Hom_{H({\bf t})}(\Lambda_2({\bf t}), \Lambda_1({\bf t})) & 
\rightarrow
& \Hom({\bf t}).
\end{array}$$

\noindent Here, the horizontal arrows are the inclusions defined above and 
the vertical arrows are the surjective restriction maps.

Notice that the space $\Hom({\bf t})$ can be represented as

$$\otimes_{i\in V_{ord}}\End_{GL(d_i)}(E_i^{\otimes p_i})\otimes
\otimes_{q\in\Omega}\End_{GL(d_q)}(E_q^{\otimes p_q}).
$$

\noindent Here $E_q\cong E_{i_q}\cong E_{j_q}$.
By Theorem 1.2 $\Hom({\bf N})$ is isomorphic to
$\otimes_{i\in V_{ord}}
K[S_{p_i}]\otimes\otimes_{q\in\Omega}K[S_{p_q}]$ since we assumed that
$N\geq r$. The kernel of the epimorphism $\Hom({\bf N})\rightarrow
\Hom({\bf t})$ is isomorphic to

$$I_{{\bf t}+1}= \sum_{i\in V_{ord}, p_i>d_i}\ldots\otimes\underbrace
{I_{d_i+1}}_{\mbox{the place of} \ K[S_{p_i}]}\otimes\ldots
+$$

$$
+\sum_{q\in\Omega, p_q>d_q}\ldots\otimes\underbrace{I_{d_q+1}}_{\mbox{the place
of} \ K[S_{p_q}]}\otimes \ldots .$$

Using Lemma 2.3 we get

\begin{pr}
If $N'\geq N\geq r$ then the epimorphism $\Hom({\bf N'})\rightarrow
\Hom({\bf N})$ is an isomorphism. The same is valid for all
epimorphisms

$$\Hom_{H({\bf N'})}(\Lambda_2({\bf N'}), \Lambda_1({\bf N'}))
\rightarrow
\Hom_{H({\bf N})}(\Lambda_2({\bf N}), \Lambda_1({\bf N})).$$
In particular, the $\bar r$-homogeneous component of the algebra
$J(Q,{\bf N})$ does not depend on the number $N$
and can be identified with the $\bar r$-homogeneous component of the free
invariant algebra $J(Q)$.
\end{pr}

Using the commutative diagram (4) from Theorem 2.1 and repeating again
the proof of Proposition 1 from \cite{zub1} we get

\begin{pr}
The $\bar{r}$-homogeneous component of the ideal $T(Q,{\bf t})$ is generated
as a vector space by
the elements $c(\phi)$, where 
$\phi \in\Hom_{H({\bf N})}(\Lambda_2({\Theta, \bf N}),
\Lambda_1({\Theta, \bf N}))\neq 0$
and $\Theta$ runs over all superpartitions of multidegree $\bar{r}$. 
In addition,
one has to require that either at least one of the modules 
$d(\Lambda_2(\Theta, {\bf N}))=\Lambda_2(\Theta, {\bf t})$, 
$d(\Lambda_1(\Theta, {\bf N}))=
\Lambda_1(\Theta, {\bf t})$ is equal to zero or $\phi\mid_{\Lambda_2(\Theta, 
{\bf t})}=0$.
\end{pr}

\begin{cor}
The algebra $J(Q,{\bf t})$ is generated by all $c(\phi)$ without any
restrictions on $\phi\in\Hom_{H({\bf t})}(\Lambda_2({\bf t}), 
\Lambda_1({\bf t}))$ or
by all $p_{{\bf N}, {\bf t}}(c(\phi')),  \phi'\in
\Hom_{H({\bf N})}(\Lambda_2({\bf N}),  \Lambda_1({\bf N}))$.
\end{cor}

\subsection{Multilinear invariants}

Denote the image of a $\phi\in\Hom({\bf t})$ in $J(\hat{Q},{\bf t})$
by $tr^*(\phi)$
and its specialization under $f$ by $tr^*(\phi, f)$.

Given operator $\sigma$ from $\Hom({\bf N})\subseteq\End_{GL(N)}
(E^{\otimes r})=K[S_r]$, where $E$ is a $N$-dimensional space,
it can be written as

\begin{equation}
\begin{array}{c}
\sum_{1\leq j_1,\ldots , j_r\leq N} (\otimes_{1\leq k\leq t}
(e_{j_k}^{t(k)})^*)\otimes (\otimes_{t+s+1\leq k\leq r}(e_{j_{k-s}}^{t(k)})^*)
\otimes (\otimes_{t+s+1\leq k\leq r}(e_{j_{k}}^{i(k)})^*)\otimes \\
\\
(\otimes_{1\leq k\leq t}e_{j_{\sigma^{-1}(k)}}^{i(k)})
\otimes(\otimes_{t+1\leq k\leq t+s} e_{j_{\sigma^{-1}(k)}}^{t(k)})\otimes
(\otimes_{t+1\leq k\leq t+s}e_{j_{\sigma^{-1}(k+s)}}^{i(k)}).
\end{array}
\end{equation}

Rearranging the tensor multipliers we see that $tr^*(\sigma)$ is equal to

$$\sum_{1\leq j_1,\ldots , j_r\leq N}(\prod_{1\leq k\leq t}y(k)_{j_k,
j_{\sigma^{-1}(k)}})(\prod_{t+1\leq k\leq t+s}y(k)_{j_{\sigma^{-1}(k)},
j_{\sigma^{-1}(k+s)}})(\prod_{t+s+1\leq k\leq r}y(k)_{j_{k-s},
j_{k}}).$$

In order to contract this sum into a product of ordinary traces
one can use the following
rule (see \cite{zub5}). We consider the {\it formal} product of pairs:

$$\prod_{1\leq k\leq t}(k, \sigma^{-1}(k))\prod_{t+1\leq k\leq t+s}
(\sigma^{-1}(k), \sigma^{-1}(k+s))\prod_{t+s+1\leq k\leq r}
(k-s, k).$$

The next step is to
partition this product into {\it cyclic} subproducts. By definition, these
subproducts are
$\prod_{1\leq f\leq l}(a_f, b_f)$ such that $b_f=a_{f+1},
1\leq f\leq l-1$ and $b_l=a_1$. If it is necessary, one can change the initial
order  of coordinates of any pair. This partition is possible because of
the
following fact: each symbol $k$ appears two times in the
original product.
Finally, each subproduct $\prod_{1\leq f\leq l}(a_f, b_f)$
corresponds to a trace $tr(Z(j_1)\ldots Z(j_l))$, where $j_f$ is the number of
the pair $(a_f, b_f)$ in the original product and $Z(j_f)$ coincides with
$Y(j_f)$ iff the initial order of the coordinates of
this pair was not changed, otherwise $Z_(j_f)=Y(j_f)^t, 1\leq f\leq l$. 

It is
more convenient for further computations to denote $Y^t$ by
$\overline{Y}$.

\begin{ex} Let $t=3,s=2,r=7,\sigma=(1726)(354)\in S_7$. The formal product
of pairs corresponding to $\sigma$ is $(16)(27)(34)(52)(31)(46)(57)$.
Decomposing it into cyclic subproducts we get $(16)(64)(43)(31)\cdot(27)
(75)(52)$. Therefore,

$$tr^*(\sigma)=tr(Y(1)\overline{Y(6)} \ \overline{Y(3)}Y(5))tr(Y(2)\overline{
Y(7)}Y(4))$$

or
 
$$tr^*(\sigma)=tr(Z(1)Z(\bar{6})Z(\bar{3})Z(5))tr(Z(2)Z(\bar{7})Z(4))$$

with respect to the notations from the introduction.
\end{ex}

Notice that if $s=0$ then $tr^*(\sigma)=tr(\sigma)$, where $tr(\sigma)=
tr(Y(a)\ldots Y(b))\ldots tr(Y(c) \\
\ldots Y(d))$ and $(a\ldots b)\ldots
(c\ldots d)$ is a cyclic decomposition of $\sigma^{-1}$.

We set

$$T(i_q)=\{j\in\hat{A}\mid t(j)=i_q\}, I(i_q)=\{j\in\hat{A}\mid i(j)=i_q\},$$
$$T(j_q)=\{j\in\hat{A}\mid i(j)=j_q\}, I(j_q)=\{j\in\hat{A}\mid t(j)=j_q\}.$$
It is clear that $\forall q\in\Omega, T(q)=T(i_q)\bigcup T(j_q), I(q)=
I(i_q)\bigcup I(j_q)$.

\begin{lm}
Any $\sigma\in S_r$ lies in $\Hom({\bf N})$ iff the following
equations are satisfied:

\begin{enumerate}
\item $\forall i\in V_{ord}, \sigma((T(i)\bigcap\hat{A}_1)\bigsqcup (T(i)
\bigcap\hat{A}_3 -s))=(I(i)\bigcap\hat{A}_1)\bigsqcup (I(i)\bigcap
\hat{A}_2 +s).$

\item $\forall q\in\Omega, \sigma((T(i_q)\bigcap\hat{A}_1)\bigsqcup
(T(i_q)\bigcap\hat{A}_3 -s)\bigsqcup T(j_q))=
(I(i_q)\bigcap\hat{A}_1)\bigsqcup (I(i_q)\bigcap
\hat{A}_2 +s)\bigsqcup I(j_q).$
\end{enumerate}
\end{lm}

Proof. For given $\sigma\in\Hom({\bf N})$ its record (5)
can be rewritten in a more refined way

$$\sum_{1\leq j_1,\ldots ,j_r\leq N}\otimes_{i\in V_{ord}}((\otimes_{1\leq
k\leq t, t(k)=i}(e^i_{j_k})^*)\otimes (\otimes_{t+1\leq k\leq t+s,
t(k+s)=i}(e^i_{j_k})^*)\otimes(\otimes_{1\leq k\leq t,
i(k)=i}e^i_{j_{\sigma^{-1}(k)}})$$

$$\otimes (\otimes_{t+1\leq k\leq t+s,
i(k)=i}e^i_{j_{\sigma^{-1}(k+s)}}))
\otimes_{q\in\Omega}((\otimes_{1\leq k\leq t,
t(k)=i_q}(e^{i_q}_{j_k})^*)
\otimes (\otimes_{t+1\leq k\leq t+s,
t(k+s)=i_q}(e^{i_q}_{j_k})^*)$$

$$\otimes (\otimes_{t+s+1\leq k\leq r,
i(k)=j_q}(e^{j_q}_{j_k})^*)
\otimes ((\otimes_{1\leq k\leq t, i(k)=i_q}e^{i_q}_{j_{\sigma^{-1}(k)}})
\otimes (\otimes_{t+1\leq k\leq t+s, i(k)=i_q}e^{i_q}_{j_{\sigma^{-1}(k+s)}})$$

$$\otimes (\otimes_{t+1\leq k\leq t+s,
t(k)=j_q}(e^{j_q}_{j_{\sigma^{-1}(k)}})).$$

It remains to notice that each factor $e^x_{yz}$ in every summand of this sum
must appear on the dual side of the same summand, that is like $(e^x_{yz})^*$.
This completes the proof.

For the sake of convenience denote the right hand side sets of these
equations by ${\cal I}(i), {\cal I}(q)$ and the left hand side sets, that is
the arguments of the substitution $\sigma$, by ${\cal T}(i), {\cal T}(q)$
respectively. Then they can be rewritten as $\sigma({\cal T}(u))=
{\cal I}(u)$, where $u\in V_{ord}\sqcup\Omega$. In other words, for any arrow 
$j$ the end of $\sigma^{-1}(j)$ coincides with the origin of $j$ (see 
\cite{zub5}).

\begin{lm}
Any $tr(Z(a_m)\ldots Z(a_1))$
occurs as a factor of some multilinear trace products from
$J(\hat{Q},\bar{t
})(1^r)$ iff $Z(a_m)\ldots Z(a_1)$ is
admissible.
\end{lm}

Proof. Fix some subproduct $UVW$ of any cyclic permutation of the product
$Z(a_m)\ldots Z(a_1)$ consisting of three
factors.
Without loss of generality one can
assume that $V=Y(j)$. Otherwise one can transpose the product
$Z(a_m)\ldots Z(a_1)$. The following list contains all admissible cases to
occupy both places around $V$.

\begin{enumerate}
\item If $j\in \hat{A}_1=[1,\ldots , t]$ then $U$ can be occupied by $Y(\sigma
(j))$. It happens iff $\sigma(j)\in \hat{A}_1$. Let $t(j)=i$. Then we have
either $j\in T(i)\bigcap\hat{A}_1$ or $i=i_q, j\in T(i_q)\bigcap\hat{A}_1$.
In both cases $\sigma(j)\in I(i)$, that is the product $Y(\sigma(j))Y(j)$ is
linked. The matrix $U$ can be equal to
$Y(j')$ or $\overline{Y(j')}$, where $j'\in\hat{A}_2=[t+1,\ldots , t+s]$.
The case $U=Y(j')$ takes place iff $\sigma(j)=j'+s\in\hat{A}_3$. It means
that either $j'\in I(i)\bigcap\hat{A}_2$ or $i=i_q, j'\in I(i_q)\bigcap
\hat{A}_2$ and in both cases the product $Y(j')Y(j)$ is also linked.

Finally, let $U=\overline{Y(j')}$. It means that $\sigma(j)=j'$. In this case
$i=i_q$ only and $j'\in I(j_q)\bigcap\hat{A}_2$, that is the product $\overline
{Y(j')}Y(j)$ is linked again.
As for $W$ the possibilities are the following: $Y(\sigma^{-1}(j)), Y(j')$ or
$\overline{Y(j')}$, where $j'\in\hat{A}_3$. As above it can easily be checked
that $VW$ is linked. Briefly one can describe all these ways of
occupying as
$\underbrace{(\hat{A}_2,\bar{\hat{A}_2}, \hat{A}_1)}_U\underbrace{\hat{A}_1}_V
\underbrace{(\hat{A}_3,\bar{\hat{A}_3},\hat{A}_1)}_W$.

Other cases are listed without any comments. The interested reader
can check them very easily.
\item If $j\in\hat{A}_2$ then either $U=\overline{Y(j')}, j'\in\hat{A}_1, t(j')
=i_q, t(j)=j_q$ or $U=Y(j'), \overline{Y(j')}, j'\in\hat{A}_3$. In the last
case
either $t(j)=i(j')$ or $t(j)=j_q, t(j')=i_q$. For $W$ we have the following
possibilities: $W=Y(j'), j'\in\hat{A}_1,$ or
$W=Y(j'), \overline{Y(j')}, j'\in\hat{A}_3$. 
The first possibility is described in the previous item, the second one
implies either $t(j')=i(j)$ or
$i(j')=j_q, i(j)=i_q$. 
Briefly,
$\underbrace{(\hat{A}_3,\bar{\hat{A}_3}, \bar{\hat{A}_1})}_U\underbrace
{\hat{A}_2}_V\underbrace{(\hat{A}_3,\bar{\hat{A}_3},\hat{A}_1)}_W$.
\item If $j\in\hat{A}_3$ then either $U=Y(j'), j'\in\hat{A}_1, t(j)=i(j')$ or
$U=Y(j'), \overline{Y(j')}, j'\in\hat{A}_2$. The last case is considered in the
second item up to some transposition. For $W$ we have the following
possibilities: $W=\overline{Y(j')}, j'\in\hat{A}_1,$ or $W=Y(j'), 
\overline{Y(j')}, j'\in\hat{A}_2$. The first possibility is described in 
the first item up to some transposition, the second possibility is described 
in the second item.
Briefly,
$\underbrace{(\hat{A}_2,\bar{\hat{A}_2}, \hat{A}_1)}_U\underbrace{\hat{A}_3}_V
\underbrace{(\hat{A}_2,\bar{\hat{A}_2},\bar{\hat{A}_1})}_W$.
\end{enumerate}

It is clear that in all cases listed above the products $UV, VW$ are
linked. The lemma is proved.

In other words, the conditions 1, 2 in Lemma 2.4 are equivalent to the 
conditions of admissibility in Lemma 2.5 .

A trace product $u=tr(Z(a_r)\ldots Z(a_k))\ldots tr(Z(a_m)\ldots Z(a_1))$
from $J(\hat{Q},{\bf N})(1^r)$ can be
written in many ways. Fix some standard record of each product
as follows. Any matrix $Z(a)$
is equal either to $Z(j)=Y(j)$ or to $Z(\bar{j})=\overline{Y(j)}$. Let us 
ascribe to $Z(a)$ this
number
$j$. The record of $tr(Z(a)\ldots Z(b))$ is called {\it right} if the matrix
with maximal number, say $j$, occupies the first place. Moreover,
$Z(a)=Y(j)$ otherwise one has to transpose the product $Z(a)\ldots Z(b)$.
Let us call $j$ by the {\it number associated to} 
$tr(Z(a)\ldots Z(b))$. The record of $u$ is called {\it right} iff all its
factors are right and their associated numbers increase on passing by this 
product from left to right.

\begin{pr}
The right trace products form a basis of the vector space
$J(\hat{Q},{\bf N})(1^r)$.
In particular, they span $J(\hat{Q},{\bf t})(1^r)$
for any ${\bf t}$.
\end{pr}

Proof. The first assertion has been proved in \cite{zub5}. The second one is
a trivial consequence of Corollary 2.1.

For the sake of convenience we will omit the symbol $tr$ in the record of
any multilinear invariant from
$J(\hat{Q},{\bf N})(1^r)$ if it does not lead to confusion.
We replace any matrix $Y(j)$ or its transposed $\overline{
Y(j)}$ by the number $j$ or its {\it transposed} $\bar{j}$, $1\leq j\leq r$.
For example, the invariant $tr(Y(1)\overline{Y(6)} \ \overline{Y(3)}
Y(5))tr(Y(2) \overline{Y(7)} Y(4))$ given above can be rewritten as
$(1\bar{6}\bar{3}5)(2\bar{7}4)$.

We suppose by definition that $\bar{\bar{i}}=i, i=1,\ldots , r$
and $[\bar{1},\bar{r}]=\{\bar{1},\ldots, \bar{r}\}$.

We reformulate the contracting rules mentioned above as follows.

\begin{pr}
Let $\sigma\in \Hom({\bf N})$ and $tr^*(\sigma)=(a\ldots b)\ldots 
(c\ldots d)$, where $a,\ldots ,b, \\ c,\ldots ,
d\in [1,r]\bigsqcup[\bar{1},\bar{r}]$. All we need is to define exactly what
is a right hand side neighbor of any symbol $j$ in a cyclic record
of $tr^*(\sigma)$?
If $j$ is an ordinary symbol, that is if $j\in [1,r]$, then we have

\begin{enumerate}
\item If $j\in \hat{A}_1$ then $(\ldots jk\ldots)$, where

$$k=\left\{\begin{array}{c}
\sigma^{-1}(j), \ \sigma^{-1}(j)\in \hat{A}_1, \\
\sigma^{-1}(j)+s, \ \sigma^{-1}(j)\in \hat{A}_2, \\
\overline{\sigma^{-1}(j)}, \ \sigma(j)\in \hat{A}_3 .
\end{array} \right. 
$$

\item If $j\in\hat{A}_2$ then $(\ldots jk\ldots )$, where

$$k=\left\{\begin{array}{c}
\sigma^{-1}(j+s), \ \sigma^{-1}(j+s)\in\hat{A}_1, \\
\sigma^{-1}(j+s)+s, \ \sigma^{-1}(j+s)\in\hat{A}_2, \\
\overline{\sigma^{-1}(j+s)}, \ \sigma^{-1}(j+s)\in\hat{A}_3 .
\end{array} \right. 
$$

\item If $j\in\hat{A_3}$ then $(\ldots jk\ldots )$, where

$$k=\left\{\begin{array}{c}
\overline{\sigma(j)}, \ \sigma(j)\in\hat{A}_1, \\
\sigma(j), \ \sigma(j)\in\hat{A}_2, \\
\overline{\sigma(j)-s}, \ \sigma(j)\in\hat{A}_3 .
\end{array} \right. 
$$
\end{enumerate}

If $j=\bar{l}$ then the corresponding rules are:

\begin{enumerate}
\item If $l\in \hat{A}_1$ then $(\ldots jk\ldots)$, where

$$k=\left\{\begin{array}{c}
\overline{\sigma(l)}, \ \sigma(l)\in \hat{A}_1, \\
\sigma(l), \ \sigma(l)\in \hat{A}_2, \\
\overline{\sigma(l)-s}, \ \sigma(l)\in \hat{A}_3 .
\end{array} \right. 
$$

\item If $l\in\hat{A}_2$ then $(\ldots jk\ldots )$, where

$$k=\left\{\begin{array}{c}
\sigma^{-1}(l), \ \sigma^{-1}(l)\in\hat{A}_1, \\
\sigma^{-1}(l)+s, \ \sigma^{-1}(l)\in\hat{A}_2, \\
\overline{\sigma^{-1}(l)}, \ \sigma^{-1}(l)\in\hat{A}_3 .
\end{array} \right. 
$$

\item If $l\in\hat{A_3}$ then $(\ldots jk\ldots )$, where

$$k=\left\{\begin{array}{c}
\overline{\sigma(l-s)}, \ \sigma(l-s)\in\hat{A}_1, \\
\sigma(l-s), \ \sigma(l-s)\in\hat{A}_2, \\
\overline{\sigma(l-s)-s}, \ \sigma(l-s)\in\hat{A}_3 .
\end{array} \right. 
$$
\end{enumerate}
\end{pr}

Proof. The check of these rules is obvious. For example, let $j=\bar{l}$ and
$\sigma^{-1}(l)\in \hat{A}_3$. Then the corresponding pair is 
$(\sigma^{-1}(l+s), \sigma^{-1}(l))$. If its right hand side neighbor (up to
the order) is $(k, \sigma^{-1}(k)), k\in\hat{A}_1,$ then either 
$\sigma^{-1}(l)=k\in\hat{A}_3$ or $l=k\in\hat{A}_2$. Both cases drive to 
a contradiction so this neighbor should be $(k, k-s), k\in\hat{A}_3$.
The other cases can be considered in the same way.

\subsection{$Z$-forms}

In the definition of the representation space of a quiver $Q$ 
of dimension ${\bf d}$ one can replace all spaces by free $Z$-modules of the
same ranks $d_1, \ldots , d_n$. We denote these modules by the same symbols 
$E_1,\ldots, E_n$. Then the free $Z$-module 
$R_Z(Q, {\bf d})=\prod_{a\in A}\Hom_Z(E_{i(a)},
E_{t(a)})$ can be regarded as a $Z$-form of $R(Q, {\bf d})$, that is
$R(Q, {\bf d})=K\otimes_Z R_Z(Q, {\bf d})$ and the dimension of the space 
$R(Q, {\bf d})$ coincides with the rank of the free $Z$-module 
$R_Z(Q, {\bf d})$. The same is true for $R_Z(Q, {\bf t})=\prod_{a\in A}
\Hom_Z(W_{i(a)},W_{t(a)})$, where as above $W_i=E_i$ iff $t_i=d_i$, otherwise
$W_i=E_i^*$.

It is clear that the ring $Z[R(Q, {\bf t})]=
Z[y_{ij}(a)\mid 1\leq j\leq d_{i(a)}, 1\leq i\leq d_{t(a)}, a\in A]$ is a
$Z$-form of $K[R(Q, {\bf t})]$. Moreover, $Z[R(Q, {\bf t})]$ can be identified
with 

$$(\prod_{a\in A_1}\otimes S(E^*_{t(a)}\otimes E_{i(a)}))
\otimes (\prod_{a\in A_2}\otimes S(E_{t(a)}\otimes
E_{i(a)}))\otimes
(\prod_{a\in A_3}\otimes S(E^*_{t(a)}\otimes E^*_{i(a)}))
$$

\noindent by the same rule as in Section 2. By Remark 1.1 any homogeneous 
component
$Z[R(Q, {\bf t})](\bar{r})$ has an ABW-filtration 

\begin{equation}
\ldots\subseteq M_{\Theta, Z}({\bf t})=M_{\Theta, Z}\subseteq\ldots
\end{equation}

\noindent such that (1) can be obtained from (6) by base change. Analogously,
for any superpartition $\Theta$ of degree $\bar{r}$ we have a homomorphism 

\begin{equation}
d_{\Theta, Z} : \Hom_Z(\Lambda_2(\Theta, Z)
, \Lambda_1(\Theta, Z))\rightarrow Z[R(Q, {\bf t})](\bar{r})
\end{equation}

\noindent such that $d_{\Theta}=K\otimes_Z d_{\Theta, Z}$. 
Here $\Lambda_1(\Theta, Z)
, \Lambda_2(\Theta, Z)$ or, more precisely, $\Lambda_1(\Theta, {\bf t}, Z),\\
\Lambda_2(\Theta, {\bf t}, Z)$
are obvious $Z$-forms of $\Lambda_1(\Theta, {\bf t})
, \Lambda_2(\Theta, {\bf t})$. If $\phi\in \Hom_Z(\Lambda_2(\Theta, Z)
, \Lambda_1(\Theta, Z))$ we denote the element $d_{\Theta, Z}(\phi)$ by
$c_Z(\phi)$. Now it is easy to guess what $tr^*_Z$ means.

We have $c(K\otimes_Z -)=K\otimes_Z c_Z(-)$, that is 
$c(K\otimes_Z \phi)=K\otimes_Z c_Z(\phi), 
\phi\in \Hom_Z(\Lambda_2(\Theta, Z), \Lambda_1(\Theta, Z))$.  
Analogously, $tr^*(K\otimes_Z -)=K\otimes_Z tr^*_Z(-)$.

Finally, as above one can define an inclusion 
$\Phi_{\Theta, Z} : \Hom_Z(\Lambda_2(\Theta, Z), 
\Lambda_1(\Theta, Z))\rightarrow B_Z({\bf t})$, where

$$B_Z({\bf t})=\Hom_Z((\otimes_{a\in A_1}
E_{t(a)}^{\otimes r_a})\otimes
(\otimes_{a\in A_3} E_{t(a)}^{\otimes r_a})\otimes (\otimes_{a\in A_3}
E_{i(a)}^{\otimes r_a}),
$$

$$(\otimes_{a\in A_1}(E_{i(a)}^{\otimes r_a})\otimes
(\otimes_{a\in A_2} E_{t(a)}^{\otimes r_a})\otimes (\otimes_{a\in A_2}
E_{i(a)}^{\otimes r_a}))
$$

\noindent such that the restriction of $K\otimes_Z \Phi_{\Theta, Z}$ on
$\Hom_{H({\bf t})}(\Lambda_2(\Theta), \Lambda_1(\Theta))$ coincides with
$\Phi_{\Theta}$.

\section{Proof of Theorem 1}

The proof of Theorem 1 is organized as follows. Using the above 
identification of the space $\Hom({\bf N})$ with a group algebra of a
product of several symmetric groups one can find some filtration of
this algebra such that the ideal $I_{{\bf t}+1}$ corresponds to a segment of 
this filtration.
We consider the above group algebra as a weight subspace of some
tensor product of symmetric algebras with respect to a torus action. 
The last tensor product has an ABW-filtration. The intersection of its terms
with our group algebra gives the required filtration of $\Hom({\bf N})$.
We notice (see Lemma 3.3 below) that the members of the last filtration are
invariant subspaces of terms of the initial ABW-filtration with respect to
an action of a reductive group which is a product of general linear groups.
By Theorem 1.1(1) a basis of $\Hom({\bf N})$ or $I_{{\bf t}+1}$ can be produced
as a union of bases of invariant subspaces of some sequential quotients of this
ABW-filtration. Following this way we get the generators of $J(Q)$ as a vector
space such that some subset of them generate the ideal $T(Q, {\bf t})$. 
It remains to simplify these generators. We reduce this 
problem to the computation of invariants of ordinary representations of 
quivers and refer to \cite{don2, zub4}.  

\subsection{Suitable generators}

Fix some $\sigma_0\in \Hom({\bf N})$. By Lemma 2.4 we have 
$\Hom({\bf N})=\sigma_0\cdot
(\otimes_{i\in V_{ord}} K[S_{{\cal T}_i}])\otimes (\otimes_{q\in\Omega}
K[S_{{\cal T}_q}])$. Moreover, the ideal
$I_{{\bf t}+1}$ is equal to

$$\sigma_0\cdot (\sum_{i\in V_{ord}, p_i>d_i}\ldots\otimes
\underbrace{I_{d_i+1}}_
{\mbox{the place of} \ K[S_{{\cal T}_i}]}\otimes\ldots
+\sum_{q\in\Omega, p_q>d_q}\ldots\otimes\underbrace{I_{d_q+1}}_{\mbox{the place
of} \ K[S_{{\cal T}_q}]}\otimes\ldots).$$

Denote by $S_{{\cal T}}$
the group $(\prod_{i\in V_{ord}} S_{{\cal T}_i})\times (\prod_{q\in\Omega}
S_{{\cal T}_q})$ and by $B({\bf t})$ the space $K\otimes_Z B_Z({\bf t})$. 

\begin{rem}
Notice that the layers of the group
$S_{\Theta_1}$ form a
a subdecomposition of the decomposition $(\bigsqcup_{i\in V_{ord}}{\cal I}(i))
\bigsqcup (\bigsqcup_{q\in\Omega}{\cal I}(q))$ and the layers of
the group $S_{\Theta_2}$ form a
subdecomposition of the decomposition $(\bigsqcup_{i\in V_{ord}}{\cal T}(i))
\bigsqcup (\bigsqcup_{q\in\Omega}{\cal T}(q))$. 
\end{rem}

\begin{lm}
The image of the module $\Hom_Z(\Lambda_2({\bf N}, Z), \Lambda_1({\bf N}, Z))$ 
in $B({\bf N})$ equals 
$\{\phi\in B({\bf N})\mid\forall \tau_1\in S_{\Theta_1}, \forall \tau_2\in 
S_{\Theta_2}, \tau_1\phi\tau_2=(-1)^{\tau_1}(-1)^{\tau_1}\phi\}$. 
The image of 
$\Hom_{H({\bf N})}(\Lambda_2({\bf N}), \Lambda_1({\bf N}))$ in the
space $\Hom(\bf N)$ is equal to
$N_{\Theta}=\{\phi\in \sigma_0\cdot K[S_{{\cal T}}]\mid\forall \tau_1\in 
S_{\Theta_1}, \forall \tau_2\in S_{\Theta_2}, \tau_1\phi\tau_2
=(-1)^{\tau_1}(-1)^{\tau_1}\phi\}$
or to
$\sigma_0\cdot\{\phi\in K[S_{{\cal T}}]\mid \forall \tau_1\in
S_{\Theta_1}, \forall \tau_2\in S_{\Theta_2}, \sigma_0^{-1}\tau_1\sigma_0
\phi\tau_2=(-1)^{\tau_1}(-1)^{\tau_1}\phi\}.$
\end{lm}

Proof. The first statement is a consequence of Lemma 1.1. We sketch the proof
of the second one and refer to \cite{zub1, zub4} for details.
It is clear that the image mentioned above is contained in $N_{\Theta}$.
On the other hand, both $\Hom_K(\Lambda_2({\bf N}), \Lambda_1({\bf N}))$ and
$B({\bf N})$ are $H({\bf N})$-modules with GF and their formal
characters do not depend on the characteristic of the ground field.
In particular, the dimensions of the spaces 
$\Hom_{H({\bf N})}(\Lambda_2({\bf N}), \Lambda_1({\bf N}))$ and
$\Hom({\bf N})$ are equal to multiplicities of the trivial character and also 
do not depend on the characteristic. It remains to notice that in the
characteristic zero case our statement is obviously true.
 
Up to the beginning of Proposition 3.2 we denote by $E_r$ a vector space of 
dimension 
$r$ with a fixed basis $e_1,\ldots , e_r$.
For any subset $T\subseteq\{1,\ldots , r\}$ denote by $E_T$ the subspace of
$E_r$ generated by all vectors $e_j,j\in T$.
Let us identify the group algebra $K[S_r]$ with a subspace of the
homogeneous
component $S^r(E_r\otimes E_r)$ by the rule $\sigma\longleftrightarrow
\prod_{i=1}^{i=r} e_{\sigma(i)}\otimes e_i$. Denote by $GL({\cal T})$ the group
$\prod_{i\in V_{ord}} GL(E_{{\cal T}(i)})\times\prod_{q\in\Omega}(E_{{\cal T}(q)})$.

We consider the space $S^r(E_r\otimes E_r)$ as a $GL(r)\times GL(r)$-module.
The group $S_r$ acts on the space $E_r$ by the rule $\sigma(e_i)=
e_{\sigma(i)}, \sigma\in S_r, 1\leq i\leq r$. In other words, we identify
the group $S_r$ with a subgroup of the group of permutation matrices by the 
rule
$\sigma\mapsto
\sum_{1\leq i\leq r} e_{\sigma(i), i}$, where $e_{kl}$ is a matrix unit which
has zero entries outside of $k$-th row or $l$-th column
but the remining entry is 1. Denote the matrix
$\sum_{1\leq i\leq r} e_{\sigma(i), i}$ by the same symbol $\sigma$.

The inclusion $K[S_r]\rightarrow S^r(E_r\otimes E_r)$ is a morphism of
$S_r\times S_r$-modules. It can easily be checked that $K[S_r]$
coincides with the weight subspace $S^r(E_r\otimes E_r)^{(1^r)\times (1^r)}$
under the induced action of the standard torus $T(r)\times T(r)$.
In the same way the space $K[S_{{\cal T}}]$ coincides with the subspace

$$((\otimes_{i\in V_{ord}} S^{p_i}(E_{{\cal T}_i}\otimes E_{{\cal T}_i}))
\otimes (\otimes_{q\in\Omega} S^{p_q}(E_{{\cal T}_q}\otimes E_{{\cal T}_q})))
^{(1^r)\times (1^r)}$$

Let $GL(\Theta_1)$ (respectively, $GL(\Theta_2)$) be a subgroup of the group 
$GL(r)$ consisting of all block
diagonal matrices  which satisfy the following requirement: if we
decompose the interval $[1, r]$ into sequential subintervals whose
lengths equal to the sizes of their blocks considered from top to bottom 
then we get the layers of the superpartition $\Theta_1$
(respectively, the layers of the superpartition
$\Theta_2$).

\begin{lm} The space $\sigma_0^{-1}\cdot N_{\Theta}$ can be identified with
$$\{g\in (\otimes_{i\in V_{ord}} S^{p_i}(E_{{\cal T}_i}\otimes E_{{\cal T}_i})
)\otimes(\otimes_{q\in\Omega} S^{p_q}(E_{{\cal T}_q}\otimes E_{{\cal T}_q}))
\mid\forall x\in GL(\Theta_1),$$
$$\forall y\in GL(\Theta_2), g^{(\sigma_0^{-1}x
\sigma_0, y)}=\det(x)\det(y)g\}.$$
\end{lm}

Proof. One has to check it on elements from $T(r)\times T(r)$ and transvections
from $GL(\Theta_1)$ and
$GL(\Theta_2)$.

Let us construct some filtration in $K[S_{{\cal T}}]$. We divide each
${\cal T}(z)$, $z\in V_{ord}\sqcup\Omega$,  into some {\it sublayers} in a 
{\it monotonic} way. In other
words, let
${\cal T}(z)=\sqcup _{1\leq j\leq l_z}\bar \beta _{zj}$, where $\max\bar \beta
_{zj_1}<\min\bar \beta _{zj_2}$ as soon as $j_1<j_2$, and $\max (\min)
\bar\beta_{zj}$ means the maximal (minimal) number from this sublayer.
Joining over all indices $z$ we obtain a decomposition of the segment
$[1, r]$.

Denote by $S_{\bar\beta }$
the Young subgroup $\prod _{i\in V_{ord}}(\prod _{1\leq j\leq
l_i}S_{\bar\beta_{ij}})\times \prod_{q\in\Omega}(\prod _{1\leq j\leq l_q}
S_{\bar\beta_{qj}})$.
As in \cite{zub1} we call this subgroup by
{\it base} subgroup.

Denote by $\Lambda^{\bar{\beta}}$ the space
$\otimes_{i\in V_{ord}, q\in\Omega}
(\otimes_{1\leq j\leq l_i} \Lambda^{p_{ij}}(E_{{\cal T}(i)}))\otimes
(\otimes_{1\leq j\leq l_q} \Lambda^{p_{qj}}(E_{{\cal T}(q)})),$
where $p_{zj}=\mid\bar{\beta}_{zj}\mid$.
The restriction of the pairing map $\delta_{\bar{\beta}}$ on the space
$\Lambda^{\bar{\beta}}\otimes\Lambda^{\bar{\beta}}$ is denoted by the same
symbol.

Repeating all arguments concerning ABW-filtrations from Section 2
we define a filtration $\{M_{\bar{\beta}}\}$
of the space
$(\otimes_{i\in V_{ord}} S^{p_i}(E_{{\cal T}_i}\otimes E_{{\cal T}_i}))\otimes
(\otimes_{q\in\Omega} S^{p_q}(E_{{\cal T}_q}\otimes E_{{\cal T}_q}))$.
For any $\bar{\beta}$ we have an exact sequence of
$GL({\cal T})\times
GL({\cal T})$-modules

\begin{equation}
0\rightarrow \ker\delta_{\bar{\beta}}\rightarrow \Lambda^{\bar{\beta}}\otimes
\Lambda^{\bar{\beta}}\rightarrow M_{\bar{\beta}}/\dot{M}_{\bar{\beta}}
\rightarrow 0.
\end{equation}

All these $GL({\cal T})\times GL({\cal T})$-modules have GF.
Denote by $G$ the group $\sigma_0^{-1}GL(\Theta_1)\sigma_0
\times GL(\Theta_2)$.

We have the filtration $\{M_{\bar\beta }^{(1^r)\times (1^r)}\}$ of the
space $K[S_{{\cal T}}]$ and $\sigma_0^{-1}\cdot I_{{\bf t}+1}$ is a 
union of members of this filtration whose indices $\bar\beta$ satisfy the
following condition: there is some $i\in V_{ord}$ or $q\in\Omega$ such that
at least one subset $\bar\beta _{ij}$ or $\bar\beta _{qj}$ has the 
cardinality $p_{ij}\geq d_i+1$ or $p_{qj}\geq d_q+1$ respectively.
Combining with Lemma 3.2 we get

\begin{lm}
The space $\sigma_0^{-1}(N_{\Theta}\bigcap
I_{{\bf t}+1})$ has the filtration
$\{(M_{\bar{\beta}}\otimes D)^G\}$,
where $D=\det^{-1}\otimes\det^{-1}$ and $\bar\beta$ satisfies the conditions
formulated above.
\end{lm}

By Remark 3.1 both groups
$\sigma_0^{-1} GL(\Theta_1)\sigma_0$ and
$GL(\Theta_2)$ are Levi subgroups of the
$GL({\cal T})$. By Theorem 1.1(4) all modules of the exact sequence 
(8) are $G$-modules with GF. Therefore, we obtain the following short exact 
sequence

$$
0\rightarrow (\ker\delta_{\bar{\beta}}\otimes D)
^G \rightarrow (\Lambda^{\bar{\beta}}\otimes\Lambda^{\bar{\beta}}\otimes D)
^G \rightarrow (M_{\bar{\beta}}/\dot{M}_{\bar{\beta}}\otimes D)
^G \rightarrow 0.
$$

\noindent In particular, all we need is to find a $(1^r)\times (1^r)$-weight 
subspace of the space
$(\Lambda^{\bar{\beta}}\otimes\Lambda^{\bar{\beta}}\otimes D)^G$ which is 
equal to
$(\Lambda^{\bar{\beta}}\otimes \ (\det)^{-1}) ^{\sigma_0^{-1} GL(\Theta_1)
\sigma_0}\otimes (\Lambda^{\bar{\beta}}\otimes (\det)^{-1})^{GL(\Theta_2)}$
and then we should compute
the image of this subspace under the pairing map $\delta_{\bar{\beta}}$.
It can easily be
checked that this subspace consists of all vectors $x$
from $(\Lambda^{\bar{\beta}}\otimes\Lambda^{\bar{\beta}})^{(1^r)\times
(1^r)}$ such that $x^{(\sigma_0^{-1}\tau_1\sigma_0, \tau_2)}=
(-1)^{\tau_1}(-1)^{\tau_2}x$, for all $\tau_1\in S_{\Theta_1}, \tau_2\in 
S_{\Theta_2}$
\cite{zub1, zub4}. Denote this subspace by $V_{\bar\beta}$.

Let $\pi\in S_{[t+1, t+s]}$ and $(a\ldots b)\ldots (c\ldots d)$ be its cyclic decomposition. Denote by $\pi +s$ the element $(a+s\ldots b+s)\ldots (c+s\ldots d+c)\in S_{[t+s+1, r]}$. Analogously, any $\pi=(a\ldots b)\ldots (c\ldots d)\in
S_{[t+s+1, r]}$ has a shifted double $\pi -s=(a-s\ldots b-s)\ldots (c-s\ldots
d-s)\in S_{[t+1, t+s]}$. 

For any Young subgroup $S_{\lambda}\leq S_{[t+1, t+s]}$
denote by $S_{\lambda +s}$ the Young subgroup of $S_{[t+s+1, r]}$ consisting of
all elements $\pi +s, \pi\in S_{\lambda}$. In the same way, $S_{\lambda -s}=
\{\pi -s\mid \pi\in S_{\lambda}\}$ if $S_{\lambda}\leq S_{[t+s+1, r]}$.  

It is clear that the groups $S_{\Theta_2}$ and
$S_{\Theta_1}$ coincide with $S_{\lambda_{A_1}}\times
S_{\gamma_{A_3}-s}\times S_{\gamma_{A_3}}$ and $S_{\lambda_{A_1}}\times
S_{\mu_{A_2}}\times S_{\mu_{A_2}+s}$ respectively. Thus any element $\pi\in
S_{\Theta_2}$ can be written as the product
$\pi_1\pi_2\pi_3$, where $\pi_1\in S_{\lambda_{A_1}}, \pi_2\in S_{\gamma_{A_3}
-s}, \pi_3\in S_{\gamma_{A_3}}$. Analogously,
any element $\pi\in S_{\Theta_1}$ can be written
as the product $\pi_1\pi_2\pi_3$, where $\pi_1\in S_{\lambda_{A_1}}, \pi_2\in
S_{\mu_{A_2}},\pi_3\in S_{\mu_{A_2}+s}$.

Denote the groups $S_{\Theta}$ and $S_{[1,t]}
\times S_{[t+1,t+s]}\times S_{[t+s+1,r]}$
by $S$ and $S_0$ respectively and define two homomorphisms $\rho_1, \rho_2$
from $S_0$ into
the group $S_r$. The first homomorphism is given by
$\pi\mapsto \pi_1\pi_2(\pi_2+s)$. The second one takes any
$\pi$ to $\pi_1(\pi_3-s)\pi_3$. 

We consider the space 
$W_{\bar{\beta}}=\{x\in (\Lambda^{\bar{\beta}}\otimes\Lambda^{\bar{\beta}})^
{(1^r)\times (1^r)}\mid \forall\tau\in S, x^{(\sigma_0^{-1}\rho_1(\tau)
\sigma_0, \rho_2(\tau))}=x\}$.
It is clear that this space contains the space
$V_{\bar{\beta}}$.

Denote by $p$ the canonical projection $\otimes_{i\in V_{ord}, 
q\in\Omega}(E_{{\cal T}(i)})^{\otimes p_i})\otimes (E_{{\cal T}(q)}
^{\otimes p_q})\rightarrow\Lambda^{\bar{\beta}}$. 

The vectors
$\bar{e}_{\sigma}=p(e_{\sigma})$ form a basis of the space
$(\Lambda^{\bar\beta})
^{(1^r)}$, where
$$e_{\sigma}=\otimes_{i\in V_{ord}, q\in\Omega}(\otimes_
{j\in {\cal T}(i)}e_{\sigma(j)})\otimes (\otimes_{j\in {\cal T}(q)}
e_{\sigma(j)})$$ 

and $\sigma$ runs over $S_{{\cal T}}/S_{\bar{\beta}}$.

\begin{pr}
The space $W_{\bar{\beta}}$ has a basis consisting of all vectors

$$\sum_{\tau\in S/(\rho_1^{-1}(S_{\bar{\beta}}^{\sigma_0\sigma_1})\bigcap
\rho_2^
{-1}(S_{\bar{\beta}}^{\sigma_2})\bigcap S)}\bar{e}_{\sigma_0^{-1}\rho_1(\tau)
\sigma_0\sigma_1}\otimes\bar{e}_{\rho_2(\tau)\sigma_2},
$$

\noindent where the pairs $(\sigma_1, \sigma_2)$ range over some subset of
$S_{{\cal T}}/S_{\bar{\beta}}\times S_{{\cal T}}/S_{\bar{\beta}}$.
\end{pr}

Proof.
The space $(\Lambda^{\bar\beta})^{(1^r)}\otimes (\Lambda^{\bar
\beta})^{(1^r)}$ has a basis $\{\bar{e}_{\sigma_1}\otimes \bar{e}_{\sigma_2}
\mid
\sigma_1, \sigma_2\in S_{{\cal T}}/S_{\bar{\beta}}\}$.
This basis is decomposed into orbits under the action of the
group $S$ by the rule $(\bar{e}_{\sigma_1}\otimes \bar{e}_{\sigma_2})^{\tau}=
\bar{e}_{\overline{\sigma_0^{-1}\rho_1(\tau)\sigma_0\sigma_1}}\otimes
\bar{e}_{\overline{\rho_2(\tau)\sigma_2}}, \tau\in S, \sigma_1, \sigma_2\in
S_{{\cal T}}/S_{\bar{\beta}}$. Therefore, it equals to

$$\bigsqcup_{(\sigma_1, \sigma_2)\in Y}\{\bar{e}_{\sigma_0^{-1}\rho_1(\tau)
\sigma_0\sigma_1}\otimes\bar{e}_{\rho_2(\tau)\sigma_2}\mid \tau\in
S/(\rho_1^{-1}(S_{\bar{\beta}}^{\sigma_0\sigma_1})\bigcap\rho_2^{-1}(S_{\bar{\beta}}
^{\sigma_2})\bigcap S)\},$$

\noindent where $Y$ is a representative set of all $S$-orbits and
$S_{\bar{\beta}}
^{\pi}=\pi S_{\bar{\beta}}\pi^{-1}, \pi\in S_r$.

It is easy to see that $\forall \tau\in S,
\sigma_0^{-1}\rho_1(\tau)\sigma_0\sigma_1 S_{\bar{\beta}}=\sigma_0^{-1}\rho_1
(\bar{\tau})\sigma_0\sigma_1 S_{\bar{\beta}}, \rho_2(\tau)\sigma_2 
S_{\bar{\beta}}=
\rho_2(\bar{\tau})\sigma_2 S_{\bar{\beta}}$. In particular, for any $\tau\in S$
we have

$$\bar{e}_{\sigma_0^{-1}\rho_1(\tau)\sigma_0\sigma_1}\otimes\bar{e}_{\rho_2
(\tau)\sigma_2}=$$

$$
(-1)^{\sigma_0^{-1}\rho_1(\tau)\sigma_0\sigma_1(\sigma_0^{-1}
\rho_1(\bar{\tau})\sigma_0\sigma_1)^{-1}}(-1)^{\rho_2(\tau)\sigma_2(\rho_2
(\bar{\tau})\sigma_2)^{-1}}\bar{e}_{\sigma_0^{-1}\rho_1(\bar{\tau})\sigma_0
\sigma_1}\otimes\bar{e}_{\rho_2(\bar{\tau})\sigma_2}=$$

$$(-1)^{\rho_1(\tau)\rho_1(\bar{\tau})^{-1}}(-1)^{\rho_2(\tau)\rho_2(\bar{\tau})
^{-1}}\bar{e}_{\sigma_0^{-1}\rho_1(\bar{\tau})\sigma_0\sigma_1}\otimes\bar{e}_
{\rho_2(\bar{\tau})\sigma_2}=\bar{e}_{\sigma_0^{-1}\rho_1(\bar{\tau})\sigma_0
\sigma_1}\otimes\bar{e}_{\rho_2(\bar{\tau})\sigma_2}.$$

This completes the proof.

We call the vectors from this proposition {\it suitable} generators.
By the same arguments one can obtain
a basis of the space $V_{\bar\beta}$. We omit these
computations and refer the interested reader to \cite{zub1, zub4}.

\begin{rem}
One can suppose that $E_r$ is a free $Z$-module with the same basis $e_1,
\ldots , e_r$. In this case $\Lambda^{\bar{\beta}}\otimes
\Lambda^{\bar{\beta}}$ is a free $Z$-module and 
we obtain the
same free generators of the free $Z$-modules $V_{\bar\beta}$ and
$W_{\bar\beta}$ as above. In
particular, any free generator of $V_{\bar\beta}$ is a sum
of suitable generators with integral coefficients \cite{zub1, zub4}.
\end{rem}

\begin{lm}
The space $N_{\Theta}\bigcap I_{{\bf t}+1}$ is generated by the elements
$$h_{\sigma_1, \sigma_2}=\sum_{\tau\in S_{\bar{\beta}}}\sum_{\pi\in S/
(\rho_1^{-1} (S_{\bar{\beta}}^{\sigma_0\sigma_1})\bigcap\rho_2^{-1}
(S_{\bar{\beta}}^{\sigma_2})\bigcap S)}(-1)^{\tau}\rho_1(\pi)\sigma_0
\sigma_1\tau\sigma_2^{-1}\rho_2(\pi)^{-1}.$$
\end{lm}

Proof. Applying the map $\delta_{\bar\beta}$ to the generators from
Proposition 3.1 we obtain the elements 

$$g_{\sigma_1, \sigma_2}=\sum_{\tau\in S_{\bar{\beta}}}\sum_{\pi\in S/
(\rho_1^{-1}(S_{\bar{\beta}}^{\sigma_0\sigma_1})\bigcap\rho_2^{-1}
(S_{\bar{\beta}}^{\sigma_2})\bigcap S)}(-1)^{\tau}\sigma_0^{-1}\rho_1(\pi)
\sigma_0\sigma_1\tau\sigma_2^{-1}\rho_2(\pi)^{-1}.$$

It remains to multiply by $\sigma_0$.

\begin{pr}
Let $\phi\in \Hom_Z(\Lambda_2({\bf t}),\Lambda_1({\bf t}))$. We have
$tr^*_Z(\Phi_{\Theta, Z}(\phi), f)=\mid
S_{\Theta}\mid c_Z(\phi)$. 
\end{pr}

Proof. Let $e^{i}_1,\ldots , e^{i}_{d_i}$ be a free basis of the module $E_i$,
$i\in V$. The
dual basis of $E_i^*$ is $(e^{i}_1)^*,\ldots , (e^{i}_{d_i})^*$.
Let us decompose the interval $[1, r]$ into subintervals by the
following rule:

$$[1, r]=(\bigsqcup_{a\in A_1}[\dot{a} , a])\bigsqcup (\bigsqcup_{a\in
A_3}[\dot{a}-s, a-s])\bigsqcup (\bigsqcup_{a\in A_3}[\dot{a}, a]),
$$

\noindent where $[\dot{a}-s, a-s]$ is equal to $[\sum_{b < a}r_b-s +1,
\sum_{b\leq a}r_b-s]$.

In the same way one can decompose the interval $[1, r]$ into other
subintervals:

$$[1, r]=(\bigsqcup_{a\in A_1}[\dot{a} , a])\bigsqcup (\bigsqcup_{a\in
A_2}[\dot{a}, a])\bigsqcup (\bigsqcup_{a\in A_2}[\dot{a}+s, a+s]),
$$

\noindent where $[\dot{a}+s, a+s]$ is equal to $[\sum_{b < a}r_b+s+1, 
\sum_{b\leq a} r_b+s]$.

Let $I, J: [1, r]\rightarrow [1, \max\limits_{i\in V}d_i]$ be
two maps such that the following conditions are satisfied:

\begin{enumerate}
\item $\forall a\in A_1, I([\dot{a}, a])\subseteq [1, d_{t(a)}]$ and
$J([\dot{a}, a])\subseteq [1, d_{i(a)}];$
\item $\forall a\in A_2, J([\dot{a}, a])\subseteq [1, d_{t(a)}]$ and
$J([\dot{a}+s, a+s])\subseteq [1, d_{i(a)}];$
\item $\forall a\in A_3, I([\dot{a}-s, a-s])\subseteq [1, d_{t(a)}]$
and $I([\dot{a}, a])\subseteq [1, d_{i(a)}].$
\end{enumerate}

Suppose that the
restrictions of the maps $I$ and $J$ on all layers of the Young subgroups
$S_{\Theta_2}$ and
$S_{\Theta_1}$
respectively are injective.
Then a typical basis vector of $\Hom_Z(\Lambda_2({\bf t}), 
\Lambda_1({\bf t}))$ is $p_{\Theta_2}(e_I^*)\otimes
p_{\Theta_1}(e_J)$
, where

$$e_I^*=(\otimes_{a\in A_1}(\otimes_{l\in [\dot{a}, a]} (e_{I(l)}^{t(a)})^*))
\otimes (\otimes_{a\in A_3}(\otimes_{l\in [\dot{a}-s, a-s]}
(e_{I(l)}^{t(a)})^*))
\otimes (\otimes_{a\in A_3}(\otimes_{l\in [\dot{a}, a]}(e_{I(l)}^{i(a)})^*))
$$

\noindent and

$$e_J=(\otimes_{a\in A_1}(\otimes_{l\in [\dot{a}, a]}e_{J(l)}^{i(a)}))
\otimes (\otimes_{a\in A_2}(\otimes_{l\in [\dot{a}, a]}e_{J(l)}^{t(a)}))
\otimes (\otimes_{a\in A_2}(\otimes_{l\in [\dot{a}+s, a+s]}e_{J(l)}^{i(a)})).
$$

The element $tr^*_Z(\Phi_{\Theta, Z}(\phi), f)$
equals

$$\sum_{\sigma_1,\sigma_2\in S_{\lambda_{A_1}}}(-1)^{\sigma_1}(-1)^{\sigma_2}
\prod_{1\leq j\leq t}y(f(j))_{I(\sigma_1(j)), J(\sigma_2(j))}\times$$

$$\sum_{\sigma_1,\sigma_2\in S_{\mu_{A_2}}}(-1)^{\sigma_1}(-1)^{\sigma_2}
\prod_{t\leq j\leq t+s}y(f(j))_{J(\sigma_1(j)), J(\sigma_2(j+s))}\times$$

$$\sum_{\sigma_1,\sigma_2\in S_{\gamma_{A_3}}}(-1)^{\sigma_1}(-1)^{\sigma_2}
\prod_{t+s\leq j\leq r}y(f(j))_{I(\sigma_1(j-s)),I(\sigma_2(j))}.$$

Ordering the factors of these products with respect to their first
subindices we get

$$\mid S_{\Theta}\mid\sum_{\sigma\in
S_{\lambda_{A_1}}}(-1)^{\sigma}
\prod_{1\leq j\leq t}y(f(j))_{I(j),J(\sigma(j))}\times$$

$$\sum_{\sigma\in S_{\lambda_{A_2}}}(-1)^{\sigma}
\prod_{t\leq j\leq t+s}y(f(j))_{J(j),J(\sigma(j+s))}\times$$

$$\sum_{\sigma\in S_{\lambda_{A_3}}}(-1)^{\sigma}
\prod_{t+s\leq j\leq r}y(f(j))_{I(j-s),I(\sigma(j))}=$$

$$=\mid S_{\Theta}\mid c_Z(\phi).$$

This concludes the proof.

\begin{cor}
Each $\frac{1}{\mid S\mid}
tr^*(h_{\sigma_1 ,\sigma_2}, f)$ belongs to $Z[R(Q, {\bf N})]$. In particular,
it can be reduced modulo any $p$. 
\end{cor}

Proof. By Lemma 3.1 we see that there is $\phi\in\Hom_Z(\Lambda_2({\bf N}, Z),
\Lambda_1({\bf N}, Z))$ such that $\Phi_{\Theta, Z}(\phi)=h_{\sigma_1, 
\sigma_2}$. By Proposition 3.2 we get  $\frac{1}{\mid S\mid}
tr^*(h_{\sigma_1 ,\sigma_2})=\frac{1}{\mid S\mid}
tr^*_Z(\Phi_{\Theta, Z}(\phi), f)=c_Z(\phi)\in Z[R(Q, {\bf N})]$.

Using Remark 3.2, Lemma 3.4 and Corollary 3.1 as well as Proposition 2.2 and
Corollary 2.1 we get

\begin{pr}
The $\bar{r}$-component of $T(Q, {\bf t})$ is generated as a vector space by
all elements $\frac{1}{\mid S\mid}
tr^*(h_{\sigma_1 ,\sigma_2}, f)$ which are also called {\it suitable},
where $S_{\bar{\beta}}$ runs over all Young
subgroups of $S_{{\cal T}}$ satisfying the condition on its layers formulated
above.
If we ignore this condition we get the generators of $J(Q)(\bar{r})$ and
mapping them into $J(Q,{\bf t})(\bar{r})$ the generators of this last
homogeneous component are obtained.
\end{pr}

\subsection{Reductions to ordinary representations of quivers}

Let us denote by $R$ the permutation
$\prod_{i\in\hat{A}_2}(i \ i+s \ \overline{i+s} \ \bar{i})$
from $S_{[1,r]\bigsqcup [\bar{1},\bar{r}]}$.
For any
$\pi\in S_{[1,r]\bigsqcup [\bar{1},\bar{r}]}$ having a cyclic decomposition
$(a\ldots b)\ldots (c\ldots d)$ denote $(\bar{a}\ldots\bar{b})\ldots
(\bar{c}\ldots\bar{d})$ by $\bar{\pi}$. We have a bijection $\iota :\pi
\mapsto
\bar{\pi}^{-1}$ on $S_{[1,r]\bigsqcup [\bar{1},\bar{r}]}$.
It is clear that this bijection induces an involution on the group
$S_{[1,r]\bigsqcup [\bar{1},\bar{r}]}$. In fact, let us denote by $a$ the 
permutation $\prod_{i\in [1, r]}(i\bar{i})$. Then $a\pi a^{-1}=\bar{\pi}$ and
$\iota(\pi)=a\pi^{-1}a^{-1}$.

\begin{lm}
Let $\sigma\in S_r$ and $tr^*(\sigma)=u=(a\ldots b)\ldots (c\ldots d)$, where
$\{a,\ldots , b, \ldots, c, \\
 \ldots , d\}$ is a subset of $[1,r]\bigsqcup [\bar{1},
\bar{r}]$ having cardinality $r$. Then $R\sigma^{-1}\bar{\sigma}R=
u\bar{u}^{-1}=u\iota(u)$.
\end{lm}

Proof. It can be easily checked that for any $j\in [1,r]\bigsqcup [\bar{1},
\bar{r}]$ its right hand side neighbors in cyclic decompositions of both
$R\sigma^{-1}\bar{\sigma}R$ and
$u$ are the same. For example, let $j=\bar{l},
l\in\hat{A}_3$. Then we have the following equations:

$$R(\bar{l})=\overline{l-s}, \bar{\sigma}(\overline{l-s})=\overline{
\sigma(l-s)}$$

and finally

$$R(\overline{\sigma(l-s)})=\left\{\begin{array}{c}
\overline{\sigma(l-s)}, \ \sigma(l-s)\in\hat{A}_1, \\
\sigma(l-s),\ \sigma(l-s)\in\hat{A}_2, \\
\overline{\sigma(l-s)-s}, \ \sigma(l-s)\in\hat{A}_3,
\end{array}\right.
$$

\noindent that is the result is the same as in the contracting rules 
defining  $u$ (see Proposition 2.5). Other cases can be checked similarly.
Thus follows that any cycle of  $R\sigma^{-1}\bar{\sigma}R$ is a cycle of
$u$ or its transposed $\bar{u}^{-1}$.
This completes the proof.

\begin{lm}
Let $\sigma\in S_r$ and $R\sigma^{-1}\bar{\sigma}R=u\iota(u)$.
Suppose that the cyclic record of $u$, including trivial cycles, contains
two symbols $i, j$ belonging to the same set $\hat{A}_l$ or
$\bar{\hat{A}}_l,l=1,2,3$. Then $(ij)u\iota((ij)u)=R\sigma'^{-1}\bar{\sigma'}R$, where either
$\sigma'=(i',j')\sigma$ or $\sigma'=\sigma (i',j')$ and $i',j'$ belong to
the same $\hat{A}_f$ or $\overline{\hat{A}}_f$, $f=1,2,3$. More precisely,

$$i',j'=\left\{\begin{array}{c}
i,j, \ i, j\in \hat{A}_1, \\
\bar{i},\bar{j}, \ i, j\in \bar{\hat{A}}_1,
\end{array}\right. $$

$$i',j'=\left\{\begin{array}{c}
i,j, \ i, j\in \hat{A}_2, \\
\overline{i+s},\overline{j+s}, \ i, j\in \bar{\hat{A}}_2,
\end{array}\right. $$

$$i',j'=\left\{\begin{array}{c}
i-s,j-s, \ i, j\in \hat{A}_3, \\
\bar{i},\bar{j}, \ i, j\in \bar{\hat{A}}_3 .
\end{array}\right. $$

\noindent In particular, $\sigma$ and $\sigma'$
have different parities.
\end{lm}

Proof. Notice that a decomposition $u\iota(u)$ of $R\sigma^{-1}\bar{\sigma}R$
is not uniquely defined.
For example, interchanging any
cycle $(a\ldots b)$ from a cyclic record of $u$ with its transposed $(\bar{b}
\ldots\bar{a})$ from a record of $\iota(u)$ we get some other decomposition
$u'\iota(u')$.
Therefore, a left factor $u$
can be defined as a part of a cyclic decomposition of $R\sigma^{-1}\bar
{\sigma}R$ depending of $r$ symbols from $[1,r]\bigsqcup [\bar{1},\bar{r}]$
which does not contain any $\iota$-invariant cycles.
Let $i,j\in\bar{\hat{A}}_3$, say $i=\bar{m}, j=\bar{n}, \ m,n\in
\hat{A}_3$. We have

$$(ij)u \ \iota((ij)u)=(\bar{m}\bar{n})R\sigma^{-1}\bar{\sigma}R(mn)=
RR^{-1}(\bar{m}\bar{n})R\sigma^{-1}\times\bar{\sigma}R(mn)R^{-1}R.$$

Further, $R^{-1}(\bar{m}\bar{n})R=(\bar{m}^{R^{-1}}\bar{n}^{R^{-1}})=(mn)$.
Thus
$(ij)u \ \iota((ij)u)=R\sigma'^{-1}\bar{\sigma'}R$, where $\sigma'=\sigma
(mn)$.
All other cases can be checked in the same way.

It remains to prove that $(ij)u$ is correctly defined. Using the identity
$(ij)(iC)(jD)=(iCjD)$, where $C, D$ are some completing fragments of these cycles,
we see that the sets of symbols involved in the records
of $u$  and $(ij)u$ correspondingly are the same. So it is enough to prove
that $(ij)u$ does not contain $\iota$-invariant cycles.

Suppose that $u=(iCjD)\ldots$. We have $(ij)u=(iC)(jD)\ldots$. If $(iC)=
\iota((iC))$ then in the cycle $(iC)$ there are two sequential
symbols like $\bar{z}, z$. Then it is true for $(iCjD)$ excepting
the case $C=C_1\bar{i}$. In the last case we have $(iCjD)=(iC_1\bar{i}
jD)$. But both cases are forbidden because of $i, j$ or $\bar{i},\bar{j}$
belongs to the same set $\hat{A}_l, l=1,2,3$, (see Lemma 2.5). 
The case $u=(iC)(jD)$ is symmetrical to the
previous one. The lemma is proved.

\begin{lm}(\cite{zub5})
Let $\pi\in S_0, \sigma\in S_r$ and $tr^*(\sigma)=u$.
Then we have $u^{\pi\times \bar{\pi}}=tr^*(\rho_1(\pi)\sigma\rho_2(\pi)^{-1})$.
\end{lm}

Proof. It is enough to prove this equation for $\pi=(ij)$, where $i,j$ lie in
$\hat{A}_1,\hat{A}_2$ or $\hat{A}_3$ simultaneously.
Let $i,j\in\hat{A}_2$. Then $\rho_1(\pi)=(ij)(i+s,j+s)$ and $\rho_2(\pi)=\id$.
We have

$$R(\bar{i}\bar{j})(\overline{i+s} \ \overline{j+s})
\sigma^{-1}\bar{\sigma}(ij)(i+s,j+s)R=((\bar{i}\bar{j})(\overline{i+s} \
\overline{j+s}))^{
\tau}R\sigma^{-1}\bar{\sigma}R((ij)(i+s,j+s))^{\tau^{-1}},$$

\noindent where $\tau=(i,i+s,\overline{i+s},\bar{i})(j,j+s,\overline{j+s},
\bar{j})$.
It remains to notice that

$$((\bar{i}\bar{j})(\overline{i+s} \ \overline{j+s}))^{\tau}=(ij)(\bar{i}
\bar{j}), ((ij)(i+s,j+s))^{\tau}=(\bar{i}\bar{j})(ij).$$

The other cases can be checked in the same way.
The lemma is proved.

We define some {\it intermediate} collection of matrices $U(l), 1\leq
l\leq m$, where $m$ is equal to the number of all layers of the group
$G=\rho_1^{-1}(S_{\bar{\beta}}^{\sigma_0\sigma_1})\bigcap\rho_2^{-1}
(S_{\bar{\beta}}
^{\sigma_2})\bigcap S$. One can define the new specialization $g$ which takes
any matrix $Y(j)$ to $U(l)$ iff $j$ belongs to the $l$-th layer of the group
$G$. It is clear that there is some specialization $h$ such that $f=h\circ g$.

\begin{lm} Every $\frac{1}{\mid S\mid} tr^*(h_{\sigma_1, \sigma_2}, f)$
is obtained from 
$\frac{1}{\mid G\mid} tr^*(\sum_{\tau\in S_{\bar{\beta}}}(-1)^{\tau}\sigma_0
\sigma_1
\tau\sigma_2^{-1}, g)$ by applying $h$.
\end{lm}

Proof. Using Lemma 3.7 we see that

$$tr^*(\rho_1(\pi)\sigma_0\sigma_1\tau
\sigma_2^{-1}\rho_2(\pi)^{-1}, f)=tr^*(\sigma_0\sigma_1\tau\sigma_2^{-1}, f)$$

\noindent because of $f\circ\pi=f$. The final computations are trivial.

Lemma 3.8 shows that without loss of generality one can assume that 
$S=G\leq\rho_1^{-1}(S_{\bar{\beta}}^
{\sigma_0\sigma_1})\bigcap\rho_2^{-1}(S_{\bar{\beta}}^{\sigma_2})$
up to some {\it gluing} of matrix
variables (see \cite{zub1, zub4}). Replacing $\sigma_0$ by $\sigma_0\sigma_1$ 
one can suppose that $\sigma_1=1$. Similarly,
replacing the group $S_{\bar{\beta}}$ by
the group $S_{\bar{\beta}}^{\sigma_2}\leq S_{{\cal T}}$ and the element
$\sigma_0$ by the element $\sigma_0\sigma_2^{-1}$ one can suppose that
$\sigma_2=1$ too.

\begin{lm}
The invariant $\frac{1}{\mid G\mid}tr^*(\sum_{\tau\in S_{\bar{\beta}}}
(-1)^{\tau}\sigma\tau, g)$ is some {\it partial linearization}
(briefly PL) of
the invariant $\frac{1}{\mid\rho_1^{-1}(S_{\bar{\beta}}^{\sigma})\bigcap
\rho_2^{-1}
(S_{\bar{\beta}})\mid}tr^*(\sum_{\tau\in S_{\bar{\beta}}}(-1)^{\tau}\sigma\tau,
f')$, where the specialization $f'$ corresponds to the group
$\rho_1^{-1}(S_{\bar{\beta}}^{\sigma})\bigcap\rho_2^{-1}
(S_{\bar{\beta}})$.
\end{lm}

Proof. By definition, $S_{f'}=\rho_1^{-1}(S_{\bar{\beta}}^{\sigma})\bigcap
\rho_2^{-1}(S_{\bar{\beta}})\leq S_0$.
Consider two layers $\alpha,
\beta$ of the group $G$ which are contained in some layer of the group
$S_{f'}\bigcap S_{[t+1, t+s]}$. For the sake of simplicity assume that
these layers have numbers $m-1, m$ correspondingly.
We define the new specialization $g'$ such $g'(j)=m-1$ iff $j\in \alpha
\bigcup\beta$ otherwise $g'(j)=g(j)$.
Let $x\in S_r$ and $tr^*(x, g')=(g'(a)\ldots g'(b))\ldots (g'(c)\ldots
g'(d))$, where $\{a,\ldots , b,
c,\ldots ,d\}$ is a subset of $[1,r]\bigcup [
\bar{1},\bar{r}]$ having cardinality $r$. By definition, $g'(\bar{
j})=\overline{g'(j)}, j\in [1,r]$.

Extracting the homogeneous summands of
degrees
$\mid\alpha\mid$ and $\mid\beta\mid$ in $U(m-1)$ and $U(m)$ respectively from
$tr^*(x, g')\mid_{U(m-1)\mapsto U(m-1)+U(m)}$ we get the sum

$$\sum_{\pi\in S_{\alpha\cup\beta}/S_{\alpha}\times S_{\beta}}(g(\pi(a))
\ldots g(\pi(b)))\ldots (g(\pi(c))\ldots g(\pi(d))).$$

\noindent Using Lemma 3.7 we see that

$$(g(\pi(a))
\ldots g(\pi(b)))\ldots (g(\pi(c))\ldots g(\pi(d)))=tr^*(\rho_1(
\pi)x, g).$$

\noindent Thus our PL of the element $\frac{1}{\mid S_{g'}\mid}tr^*(
\sum_{\tau\in
S_{\bar{\beta}}}(-1)^{\tau}\sigma\tau , g')$ is equal to

$$\frac{1}{\mid S_{g'}\mid}\sum_{\tau\in S_{\bar{\beta}}}\sum_{
\pi\in S_{\alpha\cup\beta}/S_{\alpha}\times S_{\beta}}(-1)^{\tau}tr^*
(\rho_1(\pi)\sigma\tau, g).$$

Further, $\rho_1(\pi)\in S_{\bar{\beta}}^{\sigma}$, i.e. $\rho_1(\pi)=\sigma
y\sigma^{-1}, y\in S_{\bar{\beta}}$. In particular, we get

$$\sum_{\tau\in S_{\bar{\beta}}}(-1)^{\tau}tr^*(\rho_1(\pi)
\sigma\tau, g)=
\sum_{\tau\in S_{\bar{\beta}}}(-1)^{\tau}tr^*(\sigma y\tau, g)=
\sum_{\tau\in S_{\bar{\beta}}}(-1)^{\tau}tr^*(
\sigma\tau, g),$$

\noindent since $\rho_1(\pi)$ is an even element. Therefore, our PL is equal 
to the
initial invariant $\frac{1}{\mid G\mid}tr^*(\sum_{\tau\in S_{\bar{\beta}}}
(-1)^{\tau}\sigma\tau, g)$. Repeating these arguments as many times as we need
we pass from the group $G$ to the group $S_{f'}$. This completes the proof.

Summarizing we see that up to some rearrangings, gluings of matrix variables
and PL-s the generators of $J(Q)$ ($J(Q, {\bf t})$) as well as the generators
of $T(Q, {\bf t})$ are

$$c(\phi)=\frac{1}{\mid S_f\mid}tr^*(\sum_{\tau\in S_{\bar{\beta}}}(-1)^
{\tau}\sigma\tau, f),$$

\noindent where $S_f=\rho_1^{-1}(S_{\bar{\beta}}^{\sigma})\bigcap\rho_2^{-1}
(S_{\bar{\beta}})$ and $\phi=\sum_{\tau\in S_{\bar{\beta}}}(-1)^
{\tau}\sigma\tau$. As for the generators of $T(Q, {\bf t})$ one has to impose
the condition on cardinality of layers of $S_{\bar{\beta}}$ mentioned above.
Notice that if $s=0$ then these elements are the
same as the suitable generators from \cite{zub1, zub4}.

We recall some definitions from \cite{don2}. Let $S_g\leq S_r$ be a
Young subgroup corresponding to a map $g: [1,r]\rightarrow [1,m]$.
Any sequence $p=j_1\ldots j_s$ of symbols from $[1,m]$ is said to be a
{\it primitive} cycle iff there is no proper subsequence $q$ of $p$ such
that $p$ graphically coincides with $q^k=\underbrace{q\ldots q}_k, kd=s, s>k>
1$. For any $\tau=(a\ldots b)\ldots (c\ldots d)\in S_r$ we have

$$g(\tau)=(g(a)\ldots g(b))
\ldots (g(c)\ldots g(d))=
\prod_{1\leq j\leq s_1}(p_1^{k_{1j}})\ldots\prod_{1\leq j\leq s_v}
(p_l^{k_{vj}}),$$

\noindent where each $p_i$ is a primitive cycle uniquely
defined up to cyclic permutations of its symbols,
$i=1,\ldots ,v$.
Two substitutions $\mu,\pi\in S_r$ are called $S_g$-{\it equivalent} iff
there is a sequence $\mu=\tau_1,\ldots , \tau_k=\pi$ such that for any
pair $\tau_i, \tau_{i+1}, 1\leq i\leq k-1$,
either there is $x\in S_g$ such that $\tau_{i+1}=\tau_i^x$ or
for two cycles of $\tau_i$ ($\tau_{i+1}$), say $(a\ldots b), (c\ldots d)$,
we have $(g(a)\ldots g(b))=(p^f), (g(c)\ldots g(d))=(p^d)$, where $p$ is a
primitive cycle and $\tau_{i+1}=(ac)\tau_i$ (respectively -- $\tau_i=
(ac)\tau_{i+1}$).
It can easily be checked that this relation between elements of $S_r$ is
really an equivalence. Donkin calls such equivalence class by
{\it Young superclass}.

It is clear that all permutations from the same Young superclass
$D$ have the same sets of primitive cycles. We denote each of these sets
by $P_D$.

For any Young superclass $D$ a {\it formal} invariant  
$\frac{1}{\mid S_g\mid}\sum_{x\in D}(-1)^x tr(x, g)$ can be regarded as
an invariant of $m$ $N\times N$ matrices of degree $r$ or as an element 
of the corresponding free invariant algebra due our assumption $N\geq r$.

\begin{lm}(\cite{don2})
The element $\frac{1}{\mid S_g\mid}\sum_{x\in D}(-1)^x tr(x, g)$
can be written as
a sum with integer coefficients of products of the elements $\sigma_j(p)$,
where $p\in P_D$.
\end{lm}

According to our conventions the element $\frac{1}{\mid S_g\mid}
\sum_{x\in D}(-1)^x tr(x, g)$ can be represented as a sum
$\frac{1}{\mid S_g\mid}\sum_{x\in D^{-1}}(-1)^x g(x)$, where $D^{-1}=
\{x^{-1}\mid x\in D\}$. Notice that $D^{-1}$ is also an Young superclass. 

Now, everything is prepared to prove Theorem 1.
Without loss of generality one can work in $J(Q)$ or in $J(\hat{Q})$
if it is necessary.
Let us consider any suitable generator

$$z=c(\phi)=\frac{1}{\mid S_f\mid}tr^*(\sum_{\tau\in S_{\bar{\beta}}}(-1)^
{\tau}\sigma\tau, f).$$

Fix a summand $tr^*(\sigma\tau)=u$. One can interpret the element $u$ as
an ordinary invariant depending on $r$ matrix variables
$Z(j_1),\ldots , Z(j_r), j_1,\ldots, j_r\in [1,r]\bigcup [\bar{1},\bar{r}]$ or
as a permutation from $S_{\{j_1,\ldots , j_r\}}$.
It is clear that $\{j_1,\ldots , j_r\}=T_1\bigsqcup\overline{T}_2$, where
$T_1, T_2$
are two subsets of $[1,r]$ such that $T_1\bigsqcup T_2=[1,r]$. Denote by
$S'_f$ the group $S_f^{\pi}$, where $\pi=\prod_{i\in T_2}(i\bar{i})$.
It can be easily checked that $S'_f=S_{f'}$, where $f'=(f\times f\circ a)
\mid_{j_1,\ldots, j_r}$.

I claim that Young superclass corresponding to $S'_f$, say $D$, 
which contains $u$, is a subset of
$tr^*(\sigma S_{\bar{\beta}})$. Indeed, for any $v'\in S'_f$ we have
$v'=v^{\pi}, v\in S_f$ and $u^{v'}=u^{v\times\bar{v}}=tr^*(\rho_1(v)\sigma\tau
\rho_2(v)^{-1})$ by Lemma 3.7. It remains to notice that $\sigma^{-1}\rho_1
(S_f)\sigma\leq S_{\bar{\beta}}$ and $\rho_2(S_f)\leq S_{\bar{\beta}}$.
Next, any element $(ab)$, where $a,b\in f'^{-1}(j)$ and $j$ is a symbol of some
primitive cycle belonging to $u$, has form $(ij)$ or $(\bar{i}\bar{j})$,
$(ij)\in S_f$. Using Lemma 3.6 we see that $(ab)u$ is equal to $tr^*((i'j')
\sigma\tau)$ or $tr^*(\sigma\tau (i'j'))$ and $(i'j')^{\sigma^{-1}}\in
S_{\bar{\beta}}$ or $(i'j')\in S_{\bar{\beta}}$ respectively.
For example, if $(ab)=(\bar{i}\bar{j}), i,j\in\hat{A}_2$ then $(ab)u=
tr^*((i+s,j+s)\sigma\tau)$. But any layer of $S_f\bigcap S_{[t+1,t+s]}$ has
form $\sigma(\beta_{kl})\bigcap\hat{A}_2\bigcap (\sigma(\beta_{hg})\bigcap
\hat{A}_3 -s)$. Thus $(i+s, j+s)^{\sigma^{-1}}\in S_{\bar{\beta}}
\bigcap S_{[t+s+1,r]}$.
In particular, all cycles from $P_D$ are admissible.
 
Finally, the generator $z$ can be represented as a sum of elements

$$\pm\frac{1}{\mid S_f\mid}\sum_{x\in D}(-1)^x f'(x)=
\pm\frac{1}{\mid S_{f'}\mid}\sum_{x\in D^{-1}}(-1)^x tr(x, f'),$$

\noindent where $D$ runs
over all superclasses contained in $tr^*(\sigma S_{\bar{\beta}})$.
In fact, all we need is to
prove the coincidence of signs. But for any element $u'=tr^*(\sigma\tau')$
from the Young superclass of given $u=tr^*(\sigma\tau)$
we have $(-1)^{\tau'}=(-1)^{\tau}\frac{(-1)^u}{(-1)^{u'}}$ by Lemma 3.6.
This concludes the proof.

\section{Proof of Theorem 2}

As we noticed in the introduction one can define more general supermixed
representations of quivers.
A similar definition was introduced in \cite{dw3}. Briefly speaking,
they associate with any generalized quiver of $O(n)$ ($Sp(n)$) orthogonal
(symplectic) representations of so-called {\it symmetric} quiver.
For example, typical components of orthogonal representations of a symmetric
quiver are

$$\Hom_K(V_1, V_2), \Hom_K(V_1, V_2^*), \Hom_K(V_1^*, V_2),
\Lambda^2(V)\subseteq\Hom_K(V^*, V),$$
$$\Lambda^2(V^*)\subseteq\Hom_K(V, V^*),$$
$$\Hom_K(V, W), \Hom_K(V^*, W), \Hom_K(W_1, W_2), \Lambda^2(W)\subseteq
\Hom_K(W, W).$$

The spaces $V, V_i, W, W_j$ are regarded as standard $GL(V), GL(V_i), O(W), 
O(W_j)$-modules respectively, $i=1,2, j=1,2$. These spaces are isotypical 
components of the
space $K^n$ with respect to the action of an abelian reductive subgroup $D$
of $O(n)$. The centralizer $R=Z_{O(n)}(D)$ is a product of the same
$GL(V), GL(V_i), O(W), O(W_j)$.

In the symplectic case one
has to replace the components $\Lambda^2(V), \Lambda^2(V^*), \Lambda^2(W)$
by $S^2(V), S^2(V^*), S^2(W)$ up to some identifications like
$A\longleftrightarrow AJ$ mentioned in the introduction. 
Moreover, in the last case the
groups $O(W), O(W_1), O(W_2)$ must be replaced by $Sp(W), Sp(W_1), Sp(W_2)$
correspondingly.

It is clear that our definition is more general than
Derksen-Weyman's one. For example, their definition does not include any
action of some orthogonal (symplectic) group on symmetric (skew-symmetric)
matrix component.

\begin{pr}
Let $H$ be an orthogonal or symplectic subgroup of the group $GL(n)$. 
The affine variety $GL(n)/H$ is isomorphic to
the affine variety $L$ consisting of all
non-degenerate symmetric matrices or skew-symmetric matrices
with zero diagonal entries
according to which case is considered: $H=O(n)$ or $H=Sp(n)$.
This isomorphism is induced by the map $g\mapsto gg^t$ or $g\mapsto
gJ_ng^t$ respectively.
\end{pr}

Proof. We refer to \cite{zub5} for this statement.

Notice that the left action of $GL(n)$ on $GL(n)/H$ induces the action of
$GL(n)$ on $L$ by the rule $x^g=gxg^t, x\in L, g\in GL(n)$.

Now everything is prepared to prove Theorem 2.
We describe the construction of $Q'$ step by step with respect to 
all the replacements which were used to get $S$ and $G$.

For example, let us consider the case when $G_q=Sp(d_q)$ acts on
some component $S_a\subseteq\Hom_K(V_i, V_j), a\in A, i(a)=i, t(a)=j,
V_i=V_j=K^{d_q}$ and
$S_a$ can be identified with the subspace of symmetric matrices by the rule
$A\mapsto AJ, A\in S_a$. With respect to this identification the group
$G_q=Sp(d_q)$ acts on $S_a$ by $A^g=gAg^t, g\in G_q$.

Repeating word by word the proof of Lemma 1.3 \cite{zub5} we have an
epimorphism $R^G\rightarrow K[S]^G$, where $R=K[S'\times M(d_q)]$, 
$S'$ is a product
of all components of $S$ except $S_a$, and $G_q$ acts on $M(d_q)$ by the same
rule $A^g=gAg^t$.

In fact, $S$ is a closed $G$-subvariety
of $S'\times M(d_q)$. Moreover, it is a complete intersection defined by the
relations $x_{ij}-x_{ji}=0, 1\leq i<j\leq d_q$, where $X=(x_{ij})$ is the
general matrix corresponding to the factor $M(d_q)$.

The ideal $I$ of $S$ is generated by
$G$-invariant subspace $E=\oplus_{1\leq i<j\leq d_q}K\cdot z_{ij}$,
where $z_{ij}=x_{ij}-x_{ji}, 1\leq i,j\leq d_q$.
The algebra $S(E)$ is a $G_q$-module with GF with
respect to the induced action $Z\mapsto g^{-1}Z(g^t)^{-1}, Z=(x_{ij}-
x_{ji})$. It follows immediately from  \cite{kur1, kur2}.
Using Proposition 1.3b from \cite{don7} we obtain
that $I$ is a $G$-module with GF. In particular, we have the exact sequence

$$0\rightarrow I^G\rightarrow R^G\rightarrow K[S]^G\rightarrow 0.$$

Using Proposition 4.1 we replace the group
$G_q=Sp(d_q)$ by $GL(d_q)$.
In other words, we have to add to the variety $S'\times M(d_q)$ the new
factor $GL(d_q)/Sp(d_q)$. It can be identified with a closed subvariety of
$M(d_q)^2$ consisting of all pairs of matrices $(x,y)$ such that $xy=I_{d_q}$
and both $x$ and $y$ are skew-symmetric. This subvariety is a complete
intersection again so one can use the same Proposition 1.3b from \cite{don7}.
This step was explained in
\cite{zub5} and we omit all details but briefly describe what
we get in this case.

The algebra $R^G$ is an epimorphic image of the algebra $R'^{G'}$, where
$R'=K[S'\times M(d_q)
\times M(d_q)^2]$, $G'=\times_{f, f\neq q}G_f\times GL(d_q)$ and
$GL(d_q)$ acts on the additional factor $M(d_q)^2$ by the rule
$(x,y)^g=(gxg^t, (g^t)^{-1}yg^{-1}), x, y\in M(d_q), g\in GL(d_q)$.

It means that we add to our quiver $Q$ one vertex, say with the number
$n+1$, and two arrows $b,c$ such that $i(b)=t(c)=i, t(b)=i(c)=n+1$.
Moreover, the vertex $n+1$ is occupied by the space $E_{n+1}=K^{d_q}=V$
as well as the vertex $i$ is occupied by $V^*$.
Our epimorphism is just the specialization $X(b)\mapsto J, X(c)
\mapsto -J=J^{-1}$.

As above we have the following exact sequence

$$0\rightarrow I'^{G'}\rightarrow R'^{G'}\rightarrow R^G\rightarrow 0.$$

The ideal $I'$ is generated by the $G'$-invariant subspace

$$E'=(\oplus_{1\leq i
<j\leq d_q}K\cdot z_{ij})\oplus (\oplus_{1\leq i\leq d_q}K\cdot z_i)\oplus
(\oplus_{1\leq i,j\leq d_q}K\cdot t_{ij}),$$

\noindent where $z_{ij}=x_{ij}(b)+x_{ji}(b), z_i=x_{ii}(b), 
1\leq i\neq j\leq d_q, t_{ij}=\sum_{1\leq k\leq d_q} x_{ik}(b)x_{kj}(c)-
\delta_{ij}, 1\leq i,j\leq d_q$.
All other cases can be considered in the same way as above.
This completes the proof.

\begin{rem}

There is some integer $M > 0$ depending only on ${\bf d}$ such that
whenever $char K > M$ the kernel of the epimorphism from Theorem 2 can be 
described exactly. For example, let us
consider the same case $G_q=Sp(d_q), a\in A, i(a)=i, t(a)=j, E_i=E_j=K^{d_q}$
and $S_a\subseteq\Hom_K(V_i, V_j)$ is the subspace of symmetric matrices up
to the identification $A\mapsto AJ, A\in S_a$. Using the same notations as 
above
we have the Koszul resolution

$$0\rightarrow (\Lambda^r(E)\otimes R)\rightarrow\ldots
\rightarrow (\Lambda^2(E)\otimes R)\mapsto (E\otimes R)
\rightarrow I\rightarrow 0.$$

Here $r=\dim E=\frac{d_q(d_q-1)}{2}$. Suppose $char K >r$. 
Then all $\Lambda^i(E)$
are direct summands of $E^{\otimes i}$. In particular, they are $G_q$-modules 
with GF and all members of this resolution are $G$-modules with GF.
Thus we get the exact sequence

$$0\rightarrow (\Lambda^r(E)\otimes R)^G\rightarrow\ldots
\rightarrow (\Lambda^2(E)\otimes R)^G\mapsto (E\otimes R)^G
\rightarrow I^G\rightarrow 0.$$

The space $E\otimes R$ can be considered as a homogeneous component
of an invariant ring of a supermixed representation space of some new
quiver $Q''$ having
the same set of vertices as $Q$ with the additional arrow $a''$ such that
$i(a'')=i, t(a'')=j$.
Moreover, in comparison with the previous representation
space the new one has the additional component $S_{a''}$ which is a subspace of
$M(d_q)$ consisting of all skew-symmetric matrices with respect to the same
action of $G$.

If we introduce a new general matrix $X(a'')$ corresponding to
$a''$ then $(E\otimes R)^G$ is the homogeneous component
of degree one in $X(a'')$ of the invariant algebra of this new representation 
space.
To compute the ideal $I^G$ one has to map this component to $R$ by the
rule $X(a'')\mapsto Z$. The same arguments work in the next step when we
replace $G_q=Sp(d_q)$ by $GL(d_q)$.
Notice that in the general case we have the exact sequence

$$0\rightarrow (\Delta)^G\rightarrow
(E\otimes R)^G\rightarrow I^G\mapsto H^1(G, \Delta)\rightarrow 0,$$

\noindent where $\Delta$ is the image of $\Lambda^2(E)\otimes R$ in 
$E\otimes R$.
Indeed, $H^l(G, \Delta)=0$ for all $l\geq 2$ by Lemma 1.2e (ii) from
\cite{don7}. It remains to use the long exact sequence of cohomology groups.
It was erroneously supposed in \cite{zub6} that $\Delta$ has GF and therefore
$(E\otimes R)^G\rightarrow I^G$ is an epimorphism in all characteristics.
This question is still open.
\end{rem}
\newpage
\begin{center}

Acknowledgment

\end{center}

The first variant of this article was done during author's visit to Bielefeld
University supported by DAAD. I am grateful for this support. 
Also I thank for RFFI supporting (grant N 01-01-00674).I would like to
thank Claus Ringel for his invitation and encouragement.

\end{document}